\documentclass[11pt]{amsart}

\usepackage{graphicx}
\usepackage[top=26mm, bottom=26mm, left=28mm, right=28mm]{geometry}

\newcommand{\mtx}[1]{\mathsf{#1}}

\def\p{\partial}
\def\({\left(}
\def\){\right)}
\def\[{\left[}
\def\]{\right]}

\def\xb{\boldsymbol{x}}

\def\nb{\boldsymbol{n}}

\def\yb{\boldsymbol{y}}

\newcommand{\lsp}{\vspace{2mm}}

\newtheorem{remark}{Remark}

\begin{document}

\begin{center}
\textbf{A high-order Nystr\"om discretization scheme for boundary integral
equations defined on rotationally symmetric surfaces}

\vspace{4mm}

\textit{\small P.~Young, S.~ Hao, P.G.~Martinsson}

\vspace{4mm}

\begin{minipage}{0.9\textwidth}\small
\noindent\textbf{Abstract:}
A scheme for rapidly and accurately
computing solutions to boundary integral equations (BIEs) on rotationally
symmetric surfaces in $\mathbb{R}^{3}$ is presented. The scheme uses the
Fourier transform to reduce the original BIE defined on a surface to a
sequence of BIEs defined on a generating curve for the surface. It can
handle loads that are not necessarily rotationally symmetric.
%
%
Nystr\"{o}m discretization is used to discretize the BIEs on the generating curve.
The quadrature is a high-order Gaussian rule that is modified near
the diagonal to retain high-order accuracy for singular kernels.
The reduction in dimensionality, along with the use of high-order accurate
quadratures, leads to small linear systems that can be inverted directly
via, e.g., Gaussian elimination. This makes the scheme
particularly fast in environments involving multiple right hand sides. It
is demonstrated that for BIEs associated with the Laplace and Helmholtz equations,
the kernel in the reduced equations can be evaluated very rapidly by
exploiting recursion relations for Legendre functions. Numerical
examples illustrate the performance of the scheme; in particular, it is
demonstrated that for a BIE associated with Laplace's equation on a
surface discretized using $320\,800$ points, the set-up phase of the
algorithm takes 1 minute on a standard laptop, and then solves can be
executed in 0.5 seconds.
\end{minipage}
\end{center}


\section{Introduction}
\label{sec:intro}

The premise of the paper is that it is much easier to solve a boundary integral
equation (BIE) defined on a curve in $\mathbb{R}^{2}$ than one defined on a surface
in $\mathbb{R}^{3}$. With the development of high order accurate Nystr\"om
discretization techniques \cite{Kapur:97a,alpert:99a,2001_rokhlin_kolm,2011_bremer,2011_hao_martinsson_quadrature},
it has become possible to attain
close to double precision accuracy in 2D using only a very moderate number of
degrees of freedom. This opens up the possibility of solving a BIE on
a rotationally symmetric surface with the same efficiency since such an equation can
in principle be written as a sequence of BIEs defined on a generating curve. However,
there is a practical obstacle: The kernels in the BIEs on the generating curve are
given via Fourier integrals that cannot be evaluated analytically. The principal
contribution of the present paper is to describe a set of fast methods for
constructing approximations to these kernels.

\subsection{Problem formulation}
\label{sec:intro_formulation}
\begin{figure}[htbp]
	\centering
	\includegraphics[width = 0.5\linewidth]{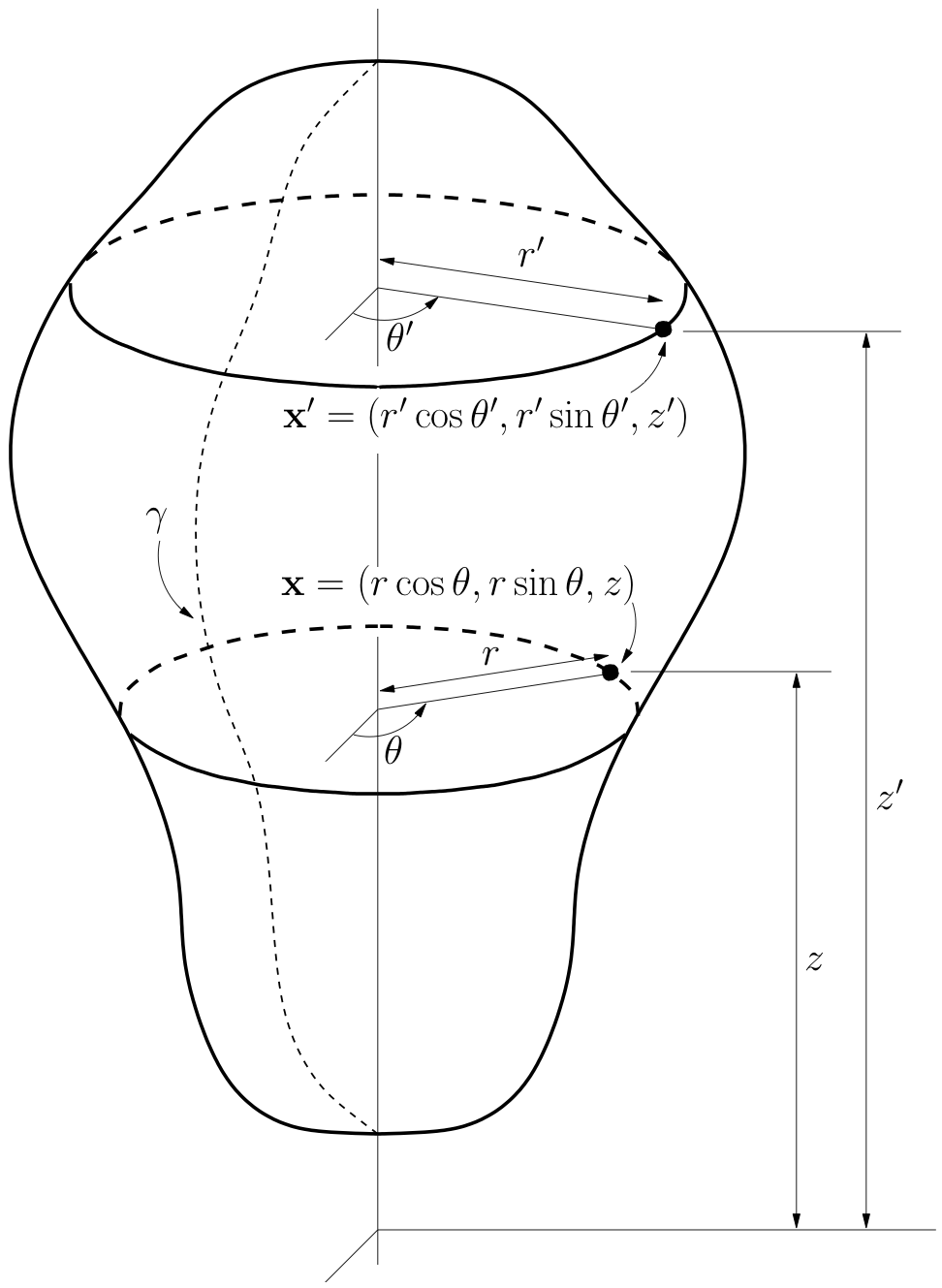}
	\caption{The axisymmetric  domain $\Gamma$ generated by the curve $\gamma$.}
	\label{fig:domain}
\end{figure}

This paper presents a numerical technique for solving boundary integral
equations (BIEs) defined on axisymmetric surfaces in $\mathbb{R}^{3}$. Specifically,
we consider second kind Fredholm equations of the form
\begin{equation}
    \sigma({\xb}) + \int_{\Gamma} k(\xb,\xb') \, \sigma(\xb') \, dA(\xb') = f(\xb), \hspace{1em} \xb \in \Gamma,
    \label{eq:fred1}
\end{equation}
under two assumptions: First, that $\Gamma$ is a surface in $\mathbb{R}^{3}$
obtained by rotating a curve $\gamma$ about an axis.
Second, that the kernel $k$ is invariant under rotation about the symmetry axis
in the sense that
\begin{equation}
    k(\xb,\xb') = k(\theta - \theta',r,z,r',z'),
    \label{eq:kerAss}
\end{equation}
where $(r,\,z,\,\theta)$ and $(r',\,z',\,\theta')$ are cylindrical
coordinates for $\xb$ and $\xb'$, respectively,
\begin{align}
\label{eq:cyl}
\xb  &= (r\,\cos\theta,\,r\,\sin\theta,\,z),\\
\xb' &= (r'\,\cos\theta',\,r'\,\sin\theta',\,z'),
\end{align}
see Figure \ref{fig:domain}. Under these assumptions, the equation (\ref{eq:fred1}), which
is defined on the two-dimensional surface $\Gamma$, can via a Fourier transform in the
azimuthal variable be recast as a sequence of equations
defined on the one-dimensional curve $\gamma$. To be precise, letting $\sigma_{n}$, $f_{n}$,
and $k_{n}$ denote the Fourier coefficients of $\sigma$, $f$, and $k$, respectively (so
that (\ref{eq:forSeries1}), (\ref{eq:forSeries2}), and (\ref{eq:forSeries3}) hold), the
equation (\ref{eq:fred1}) is equivalent to the sequence of equations
\begin{equation}
\label{eq:fred2}
    \sigma_{n}(r,z) + \sqrt{2 \pi} \int_{\gamma} k_{n}(r,z,r',z') \,\sigma_{n}(r',z') \,r'\,  dl(r',z') = f_{n}(r,z),
    \hspace{1em} (r,z) \in \gamma, \hspace{0.5em} n \in \mathbb{Z}.
\end{equation} Whenever $f$ can be represented with a moderate number of Fourier
modes, the formula (\ref{eq:fred2}) provides an efficient technique for computing
the corresponding modes of $\sigma$. The conversion of (\ref{eq:fred1}) to
(\ref{eq:fred2}) appears in, e.g., \cite{Rizzo:79a}, and is described in
detail in Section \ref{sec:sepVar}. Note that the conversion procedure does \textit{not}
require the data function $f$ to be rotationally symmetric.

\subsection{Applications and prior work}
Equations like (\ref{eq:fred1}) arise in many areas of
mathematical physics and engineering, commonly as reformulations of
elliptic partial differential equations. Advantages of a BIE
approach include a reduction in dimensionality,
often a radical improvement in the conditioning of the mathematical
equation to be solved, a natural way of handling problems defined on
exterior domains, and a relative ease in implementing high-order
discretization schemes, see, e.g., \cite{Atkinson:97a}.

The observation that BIEs on rotationally symmetric surfaces can conveniently be
solved by recasting them as a sequence of BIEs on a generating curve has previously
been exploited in the context of stress analysis \cite{Bakr:85a}, scattering
\cite{Fleming:04a,Kuijpers:97a,Soenarko:93a,Tsinopoulos:99a,Wang:97a}, and potential
theory \cite{Gupta:1979a,Provatidis:98a,Rizzo:79a,Shippy:80a}.  Most of these
approaches have relied on collocation or Galerkin discretizations and have generally
used low-order accurate discretizations.

\subsection{Kernel evaluations} \label{sec:intro_kernel} A complication of the
axisymmetric formulation is the need to determine the kernels $k_{n}$ in
(\ref{eq:fred2}). Each kernel $k_{n}$ is defined as a Fourier integral of the
original kernel function $k$ in the azimuthal variable $\theta$ in
(\ref{eq:kerAss}), cf.~(\ref{eq:def_kn}), that cannot be evaluated analytically, and
would be too expensive to approximate via standard quadrature techniques. The FFT
can be used to a accelerate the computation in certain regimes. For many points,
however, the function which is to be transformed is sharply peaked, and the FFT
would at these points yield inaccurate results.

\subsection{Principal contributions of present work} This paper resolves the
difficulty of computing the kernel functions $k_{n}$ described in Section
\ref{sec:intro_kernel}. For the case of kernels associated with the Laplace
equation, it provides analytic recursion relations that are valid precisely in the
regions where the FFT loses accuracy. The kernels associated with the Helmholtz
equation can then be obtained via a perturbation technique.

The paper also describes a high-order Nystr\"om discretization
of the BIEs (\ref{eq:fred2}) that provides far higher accuracy and speed
than previously published methods. The
discretization scheme converges fast enough that for simple generating curves,
a relative accuracy of  $10^{-10}$ is obtained using as few as a hundred points,
cf.~Section \ref{sec:num}.
The rapid convergence of the discretization leads to linear systems of
small size that can be solved \emph{directly} via, e.g.,
Gaussian elimination, making the algorithm particularly effective in
environments involving multiple right hand sides or when the linear
system is challenging for iterative solvers (as happens for many scattering problems).

Finally, the efficient techniques for evaluating the fundamental solutions to the
Laplace and Helmholtz equations in an axisymmetric environment have applications
beyond solving boundary integral equations, for details see Section
\ref{sec:fundamentalsolution}.

\subsection{Asymptotic costs}
To describe the asymptotic complexity of
the method, we let $N_{\rm tot}$ denote the total number of discretization
points, and assume that as $N_{\rm tot}$ grows, the number of Fourier
modes required to resolve the solution scales proportionate to the
number of discretization points required along the generating curve
$\gamma$.
Then the asymptotic cost of solving (\ref{eq:fred1}) for a
single right-hand side $f$ is $O(N_{\rm tot}^{2})$.
If additional
right hand sides are given, any subsequent solve requires only
$O(N_{\rm tot}^{3/2})$ operations. Numerical experiments presented in Section
\ref{sec:num} indicate that the constants of proportionality in the two
estimates are moderate. For instance, in a simulation with
$N_{\rm tot} = 320\,800$, the scheme requires 1 minute for the first
solve, and $0.49$ seconds for each additional right hand side (on a standard
laptop).

Observe that since a high-order discretization scheme is used, even complicated
geometries can be resolved to high accuracy with a moderate number $N_{\rm tot}$ of
points.

\subsection{Outline}
Section \ref{sec:sepVar} provides details on the conversion of a BIE on a rotationally
symmetric surface to a sequence of BIEs on a generating curve.
Section \ref{sec:disc_BIE} describes a high order methods for discretizing a BIE on a curve.
Section \ref{sec:genAlg} summarizes the algorithm and estimates its computational costs.
Section \ref{sec:laplace} describes how to rapidly evaluate the kernels associated with
the Laplace equation, and then Section \ref{sec:kerExtension} deals with the Helmholtz case.
Section \ref{sec:fundamentalsolution} describes other applications of the kernel evaluation techniques.
Section \ref{sec:num} illustrates the performance of the proposed method via numerical
experiments.


\section{Fourier Representation of BIE}
\label{sec:sepVar}

Consider the BIE (\ref{eq:fred1}) under the assumptions on
rotational symmetry stated in Section \ref{sec:intro_formulation}
(i.e.~$\Gamma$ is a rotationally symmetric surface generated by
a curve $\gamma$ and that $k$ is a rotationally symmetric kernel).
Cylindrical coordinates $(r,z,\theta)$ are introduced as specified
in (\ref{eq:cyl}). We write $\Gamma = \gamma \times \mathbb{T}$
where $\mathbb{T}$ is the one-dimensional torus, usually parameterized
by $\theta \in (-\pi,\,\pi]$.

\subsection{Separation of Variables}
We define for $n \in \mathbb{Z}$ the functions $f_{n}$, $\sigma_{n}$, and
$k_{n}$ via
\begin{align}
\label{eq:def_fn}
    f_{n}(r,z)            &= \int_{\mathbb{T}} \frac{e^{-in\theta}}{\sqrt{2\pi}}\, f(\theta ,r,z) \, d\theta,\\
\label{eq:def_sigman}
    \sigma_{n}(r,z)       &= \int_{\mathbb{T}} \frac{e^{-in\theta}}{\sqrt{2\pi}}\, \sigma(\theta,r,z)\, d\theta,\\
\label{eq:def_kn}
    k_{n}(r,z,r',z') &= \int_{\mathbb{T}} \frac{e^{-in\theta}}{\sqrt{2\pi}}\,k(\theta,r,z,r',z')\,d\theta.
\end{align}
Formulas (\ref{eq:def_fn}), (\ref{eq:def_sigman}), and
(\ref{eq:def_kn}) define $f_{n}$, $\sigma_{n}$, and $k_{n}$ as the
coefficients in the Fourier series of the functions $f$, $\sigma$, and
$k$ about the azimuthal variable,
\begin{align}
\label{eq:forSeries1}
f(\xb)      &= \sum_{n \in \mathbb{Z}}\frac{e^{in\theta}}{\sqrt{2\pi}}\,f_{n}(r,z),\\
\label{eq:forSeries2}
\sigma(\xb) &= \sum_{n \in \mathbb{Z}}\frac{e^{in\theta}}{\sqrt{2\pi}}\,\sigma_{n}(r,z),\\
\label{eq:forSeries3} k(\xb,\xb') = k(\theta-\theta',r,z,r',z') &=
    \sum_{n \in \mathbb{Z}}\frac{e^{in(\theta-\theta')}}{\sqrt{2\pi}}\,k_{n}(r,z,r',z').
\end{align}

To determine the Fourier representation of (\ref{eq:fred1}), we multiply
the equation by $e^{-in\theta}/\sqrt{2\pi}$ and integrate $\theta$ over
$\mathbb{T}$. Equation (\ref{eq:fred1}) can then be said
to be equivalent to the sequence of equations
\begin{equation}
\label{eq:step1} \sigma_{n}(r,z) + \int_{\gamma\times\mathbb{T}}
\left[\int_{\mathbb{T}}\frac{e^{-in\theta}}{\sqrt{2\pi}}\,k(\xb,\xb')\,d\theta\right]\,\sigma(\xb')\,dA(\xb')
= f_{n}(r,z),\qquad n \in \mathbb{Z}.
\end{equation}
Invoking (\ref{eq:forSeries3}), we evaluate the bracketed factor in (\ref{eq:step1}) as
\begin{multline}
\label{eq:step2}
\int_{\mathbb{T}}\frac{e^{-in\theta}}{\sqrt{2\pi}}\,k(\xb,\xb')\,d\theta
= \int_{\mathbb{T}}\frac{e^{-in\theta}}{\sqrt{2\pi}}\,k(\theta-\theta',r,z,r',z')\,d\theta\\
= e^{-in\theta'}\int_{\mathbb{T}}\frac{e^{-in(\theta-\theta')}}{\sqrt{2\pi}}\,k(\theta-\theta',r,z,r',z')\,d\theta
= e^{-in\theta'}\,k_{n}(r,z,r',z').
\end{multline}
Inserting (\ref{eq:step2}) into (\ref{eq:step1}) and executing the integration of $\theta'$ over $\mathbb{T}$, we find
that (\ref{eq:fred1}) is equivalent to the sequence of equations
\begin{equation}
    \label{eq:step4}
    \sigma_{n}(r,z) + \sqrt{2\pi} \int_{\gamma}\,k_{n}(r,z,r',z')\,\sigma_{n}(r',z')\,r' \,dl(r',z') = f_{n}(r,z),
     \hspace{1em} n \in \mathbb{Z}.
\end{equation}
For future reference, we define for $n\in \mathbb{Z}$ the boundary integral operators $\mathcal{K}_{n}$ via
\begin{equation}
\label{eq:def_Kn}
[\mathcal{K}_{n} \, \sigma_{n}](r,z) = \sqrt{2\pi}\int_{\gamma} k_{n}(r,z,r',z')\,\sigma_{n}(r',z')\,r'\,dl(r',z').
\end{equation}
Then equation (\ref{eq:step4}) can be written
\begin{equation}
\label{eq:step5}
\bigl(I + \mathcal{K}_{n})\,\sigma_{n} = f_{n},\qquad n\in \mathbb{Z}.
\end{equation}
When each operator $I + \mathcal{K}_{n}$ is continuously invertible,
we write the solution of (\ref{eq:fred1}) as
\begin{equation}
\label{eq:solution_op}
\sigma(r,z,\theta)  = \sum_{n\in \mathbb{Z}} \frac{e^{in\theta}}{\sqrt{2\pi}} [(I + \mathcal{K}_{n})^{-1} f_{n}](r,z).
\end{equation}



\subsection{Truncation of the Fourier series} \label{sec:trunction} When evaluating
the solution operator (\ref{eq:solution_op}) in practice, we will choose a
truncation parameter $N$, and evaluate only the lowest $2N+1$ Fourier modes. If $N$
is chosen so that the given function $f$ is well-represented by its lowest $2N+1$
Fourier modes, then in typical environments the solution obtained by truncating the
sum (\ref{eq:solution_op}) will also be accurate. To substantiate this claim,
suppose that $\varepsilon$ is a given tolerance, and that $N$ has been chosen so
that \begin{equation} \label{eq:error_in_f}
    ||f - \sum_{n=-N}^{N} \frac{e^{in\theta}}{\sqrt{2\pi}} f_{n}|| \leq \varepsilon,
\end{equation}
We define an approximate solution via
\begin{equation}
\label{eq:def_sigma_approx}
\sigma_{\textrm{approx}} = \sum_{n = -N}^{N} \frac{e^{in\theta}}{\sqrt{2\pi}} (I + \mathcal{K}_{n})^{-1} f_{n}.
\end{equation}
From Parseval's identity, we then find that the error in the solution satisfies
\begin{align*}
|| \sigma - \sigma_{\textrm{approx}} ||^{2}
= \sum_{|n| > N}||(I + \mathcal{K}_{n})^{-1} f_{n}||^{2}
\leq \sum_{|n| > N}||(I + \mathcal{K}_{n})^{-1}||^{2} \, ||f_{n}||^{2} \\
\leq \left(\max_{|n| > N}||(I + \mathcal{K}_{n})^{-1}||^{2}\right)\sum_{|n| > N} ||f_{n}||^{2}
\leq \left(\max_{|n| > N}||(I + \mathcal{K}_{n})^{-1}||^{2}\right)\varepsilon^{2}.
\end{align*}
It is typically the case that the kernel $k(\xb,\xb')$ has enough smoothness that the
Fourier modes $k_{n}(r,z,r',z')$ decay as $n\rightarrow \infty$. Then
$||\mathcal{K}_{n}|| \rightarrow 0$ as $n \rightarrow \infty$ and $||(I + \mathcal{K}_{n})^{-1}|| \rightarrow 1$.
Thus, an accurate approximation of $f$ leads to an approximation in $\sigma$ that is of the same order of accuracy.
Figure \ref{fig:cond} illustrates how fast this convergence is for Laplace's equation
(note that in the case illustrated, the original equation is $\tfrac{1}{2}\sigma + \mathcal{K}\sigma = f$,
and it is shown that $||(\tfrac{1}{2}I + \mathcal{K}_{n})^{-1}|| \rightarrow 1/2$).


\section{Nystr\"om discretization of BIEs on the generating curve}
\label{sec:disc_BIE}

We discretize the BIEs (\ref{eq:fred2}) defined on the generating
curve $\gamma$ using a Nystr\"om scheme. In describing the scheme,
we keep the formulas uncluttered by discussing a generic integral
equation
$$
\sigma(\xb) + \int_{\gamma}k(\xb,\xb')\,\sigma(\xb')\,dl(\xb') = f(\xb),\qquad \xb \in \gamma,
$$
where $\gamma$ is a simple smooth curve in the plane, and $k$ is a
weakly singular kernel function.

\subsection{Quadrature nodes}
Consider a quadrature rule on $\gamma$ with nodes $\{\xb_{i}\}_{i=1}^{I} \subset \gamma$
and weights $\{w_{i}\}_{i=1}^{I}$. In other words, for a sufficiently
smooth function $\varphi$ on $\gamma$,
\begin{equation}
\label{eq:gamma_quad}
\int_{\gamma}\varphi(\xb)\,dl(\xb) \approx
\sum_{i=1}^{I} \varphi(\xb_{i})\,w_{i}.
\end{equation}
For the experiments in this paper, we use a composite Gaussian rule
with $10$ points per panel. Such a rule admits for local refinement,
and can easily be modified to accommodate contours with corners
that are only piece-wise smooth.

\subsection{A simplistic Nystr\"om scheme}
The Nystr\"om discretization of (\ref{eq:gamma_quad}) corresponding to a
quadrature with nodes $\{\xb_{i}\}_{i=1}^{I}$ takes the form
\begin{equation}
\label{eq:nystrom}
\sigma_{i} + \sum_{j=1}^{I} a_{i,j}\,\sigma_{j} = f(\xb_{i}),
\qquad i = 1,\,2,\,3,\,\dots,\,I,
\end{equation}
where $\{a_{i,j}\}_{i,j=1}^{I}$ are coefficients such that
\begin{equation}
\label{eq:nystrom_goal}
\int_{\gamma}k(\xb_{i},\xb')\,\sigma(\xb')\,dl(\xb') \approx
\sum_{j=1}^{I} a_{i,j}\,\sigma(\xb_{j}),
\qquad i = 1,\,2,\,3,\,\dots,\,I.
\end{equation}
The solution of (\ref{eq:nystrom}) is a vector $\boldsymbol{\sigma} =
[\sigma_{i}]_{i=1}^{I}$ such that each $\sigma_{i}$ is an approximation
to $\sigma(\xb_{i})$.

A simplistic way to construct coefficients $a_{i,j}$ so that (\ref{eq:nystrom_goal}) holds
is to simply apply the rule (\ref{eq:gamma_quad}) to each function
$\xb' \mapsto k(\xb_{i},\xb')\,\sigma(\xb')$
whence
\begin{equation}
\label{eq:aij_simple}
a_{i,j} = k(\xb_{i},\xb_{j})\,w_{j}.
\end{equation}
This generally results in low accuracy since the kernel $k(\xb,\xb')$
has a singularity at the diagonal. However, the formula (\ref{eq:aij_simple})
has the great advantage that constructing each $a_{i,j}$ costs no more
than a kernel evaluation; we seek to preserve this property for as many
elements as possible.

\subsection{High-order accurate Nystr\"om discretization}
\label{sec:nystrom_highorder}
It is possible to construct a high-order discretization that preserves the
simple formula (\ref{eq:aij_simple}) for the vast majority of coefficients $a_{i,j}$
\cite{2011_hao_martinsson_quadrature}.
The only coefficients that need to be modified are those for which the target point
$\xb_{i}$ is \textit{near}\footnote{To be precise, we say that a point $\xb_{i}$ is \textit{near}
a panel $\tau$ if $\xb_{i}$ is located inside a circle that is concentric to the smallest
circle enclosing $\tau$, but of twice the radius.}
the panel $\tau$ holding $\xb_{j}$. In this case, $a_{i,j}$ is conceptually
constructed as follows: First map the pointwise values $\{\sigma_{j}\}_{\xb_{j} \in \tau}$
to their unique polynomial interpolant on $\tau$, then integrate this polynomial
against the singular (or sharply peaked) functions $\xb' \mapsto k(\xb_{i},\xb')$
using the quadratures of \cite{2001_rokhlin_kolm}
Operationally, the end result is that $a_{i,j}$ is given by
\begin{equation}
\label{eq:aij_highorder}
a_{i,j} = \left\{\begin{array}{ll}
\sum_{p=1}^{m}k(\xb_{i},\yb_{i,j,p})\,v_{i,j,p}\quad&\mbox{if }\xb_{i}\mbox{ and }\xb_{j}\mbox{ are near},\\
k(\xb_{i},\xb_{j})\,w_{j}\quad&\mbox{if }\xb_{i}\mbox{ and }\xb_{j}\mbox{ are not near}.
\end{array}\right.
\end{equation}
In (\ref{eq:aij_highorder}),
$m$ is a small integer (roughly equal to the order of the Gaussian quadrature),
the numbers $v_{i,j,p}$ are coefficients that depend on $\gamma$ but not on $k$,
and $\yb_{i,j,p}$ are points on $\gamma$ located in the same panel as $\xb_{j}$
(in fact, $\yb_{i,j,p} = \yb_{i,j',p}$ when $j$ and $j'$ belong to the same panel,
so the number of kernel evaluations required is less than it seems).
For details, see \cite{2011_hao_martinsson_quadrature}.


\section{The full algorithm}
\label{sec:genAlg}

\subsection{Overview}
\label{sec:summary}
At this point, we have shown how to convert a BIE defined on
an axisymmetric surface in $\mathbb{R}^{3}$ to a sequence of equations defined on
a curve in $\mathbb{R}^{2}$ (Section \ref{sec:sepVar}), and then how to discretize
each of these reduced equations (Section \ref{sec:disc_BIE}). Putting the components
together, we obtain the following algorithm for solving (\ref{eq:fred1}) to within
some preset tolerance $\varepsilon$:

\lsp

\begin{enumerate}
\renewcommand{\labelenumi}{\roman{enumi})}
\item Given $f$,
determine a truncation parameter $N$ such that
$||f - \sum_{n=-N}^{N}\frac{e^{in\theta}}{\sqrt{2\pi}}\,f_{n}|| \leq \varepsilon$.

\lsp

\item Fix a quadrature rule for $\gamma$ with nodes $\{(r_{i},z_{i})\}_{i=1}^{I} \subset \gamma$
and form for each Fourier mode $n = -N,\,-N \,+1,\,-N+2,\,\dots,\,N$ the corresponding Nystr\"om
discretization as described in Section \ref{sec:disc_BIE}. The number of nodes $I$ must be
picked to meet the computational tolerance $\varepsilon$. Denote the resulting coefficient
matrices $\{\mtx{A}^{(n)}\}_{n=-N}^{N}$.

\lsp

\item Evaluate via the FFT the terms $\{f_{n}(r_{i},z_{i})\}_{n=-N}^{N}$ (as defined by (\ref{eq:def_fn}))
for $i = 1,\,2,\,3,\,\dots,\,I$.

\item
Solve the equation $(\mtx{I} + \mtx{A}^{(n)})\,\sigma_{n} = f_{n}$
for $n = -N,\,-N+1,\,-N+2,\,\dots,\,N$.

\lsp

\item Construct $\sigma_{\rm approx}$
using formula (\ref{eq:def_sigma_approx}) evaluated via the FFT.

\end{enumerate}

\lsp

The construction of the matrices $\mtx{A}^{(n)}$ in Step ii can be accelerated using the FFT
(as described in Section \ref{sec:fastform_An}), but even with such acceleration, it
is typically by a wide margin the most expensive part of the algorithm. However, this
step needs to be performed only once for any given geometry.
The method therefore becomes particularly efficient when (\ref{eq:fred1}) needs to be
solved for a sequence of right-hand sides.
In this case, it may be worth the cost to pre-compute the inverse (or LU-factorization) of each matrix $\mtx{I} + \mtx{A}^{(n)}$.


\subsection{Cost of computing the coefficient matrices}
\label{sec:fastform_An}
For each of the $2N+1$ Fourier modes, we need to construct an $I\times I$ matrix $\mtx{A}^{(n)}$
with entries $a_{i,j}^{(n)}$. These entries are derived from the kernel functions $k_{n}$
defined by (\ref{eq:def_kn}). Note that whenever $(r,z) \neq (r',z')$, the function
$\theta \mapsto k(\theta,r,z,r',z')$ is $C^{\infty}$, but that as
$(r',z') \rightarrow (r,z)$ it develops a progressively sharper peak around $\theta = 0$.

For two nodes $(r_{i},z_{i})$ and $(r_{j},z_{j})$ that are ``not near'' (in the
sense defined in Section \ref{sec:nystrom_highorder}) the matrix entries are given by the
formula
\begin{equation}
\label{eq:aijn_far}
a_{i,j}^{(n)} = k_{n}(r_{i},z_{i},r_{j},z_{j})\,w_{j}
\end{equation}
where $k_{n}$ is given by (\ref{eq:def_kn}). Using the FFT, all $2N+1$
entries can be evaluated at once in $O(N \log N)$ operations. The FFT implicitly
evaluates the integrals (\ref{eq:def_kn}) via a trapezoidal rule which is highly
accurate since the points $(r_{i},z_{i})$ and $(r_{j},z_{j})$ are well-separated
on the curve $\gamma$.

For two nodes $(r_{i},z_{i})$ and $(r_{j},z_{j})$ that are not well-separated
on $\gamma$, evaluating $a_{i,j}^{(n)}$ is dicier. The first complication is
that we must now use the corrected formula, cf.~(\ref{eq:aij_highorder}),
\begin{equation}
\label{eq:aijn_close}
a_{i,j}^{(n)} =
\sum_{p=1}^{m}k_{n}(r_{i},z_{i},r_{i,j,p},z_{i,j,p})\,v_{i,j,p}.
\end{equation}
The second complication is that the FFT acceleration for computing
the kernels $\{k_{n}\}_{n=-N}^{N}$ jointly no longer works since
the integrand in (\ref{eq:def_kn}) is too peaked for the simplistic
trapezoidal rule implicit in the FFT. Fortunately, it turns out
that for the kernels we are most interested in (the single and double
layer kernels associated with the Laplace and Helmholtz equations),
the sequence $\{k_{n}\}_{n=-N}^{N}$ can be evaluated very efficiently
via certain recurrence relations as described in Sections \ref{sec:laplace}
and \ref{sec:kerExtension} (for the Laplace and Helmholtz equations,
respectively). Happily, the analytic formulas are stable precisely
in the region where the FFT becomes inaccurate.


\subsection{Computational Costs}
\label{sec:compcosts}
The asymptotic cost of the algorithm described in Section \ref{sec:summary}
will be expressed in terms of the number $N$ of Fourier modes required,
and the number $I$ of discretization points required along $\gamma$.
The total cost can be split into three components:

\lsp

\begin{enumerate}
\item \textit{Cost of forming the matrices $\{\mtx{A}^{(n)}\}_{n=-N}^{N}$:}
We need to form $2N+1$ matrices, each of size $I\times I$. For $O(I^{2})$
entries in each matrix, the formula (\ref{eq:aijn_far}) applies and using
the FFT, all $2N+1$ entries can be computed at once at cost $O(N\log N)$.
For $O(I)$ entries close to the diagonal, the formula (\ref{eq:aijn_close})
applies, and all the $2N+1$ entries can be computed at once at cost $O(N)$
using the recursion relations in Sections \ref{sec:laplace} and \ref{sec:kerExtension}.
The total cost of this step is therefore $O(I^{2}N\log N)$.

\lsp

\item \textit{Cost of transforming functions from physical space to Fourier space and back:}
The boundary data $f$ must be decomposed into Fourier modes $\{f_{n}\}_{n=-N}^{N}$, and after
the linear systems $(\mtx{I} + \mtx{A}^{(n)})\sigma_{n} = f_{n}$ have been solved, the
Fourier modes $\{\sigma_{n}\}_{n=-N}^{N}$ must be transformed back to physical space.
The asymptotic cost is $O(IN\log(N))$.

\lsp

\item \textit{Cost of solving the linear systems $(\mtx{I} + \mtx{A}^{(n)})\,\sigma_{n} = f_{n}$:}
Using a direct solver such as Gaussian elimination, the asymptotic cost is $O(I^{3}N)$.
\end{enumerate}

\lsp

We make some practical observations:
\begin{itemize}
\item The cost of executing FFTs is minimal and is dwarfed by the remaining costs.

\lsp

\item The scheme is highly efficient in situations where the same equation
needs to be solved for a sequence of different right hand sides. In this
situation, one factors the matrices $\mtx{I} + \mtx{A}^{(n)}$ once at cost
$O(I^{3}N)$, and then the cost of processing an additional right hand
side is only $O(I^{2}N + IN\log N)$ with a very small constant of proportionality.

\lsp

\item To elucidate the computational costs, let us express them in terms
of the total number of discretization points $N_{\rm tot}$ under the simplifying
assumption that $I \sim N$. Since $N_{\rm tot} = IN$, we find $I \sim N_{\rm tot}^{1/2}$
and $N \sim N_{\rm tot}^{1/2}$. Then:
\begin{center}
\begin{tabular}{lll}
Cost of setting up linear systems: & $O(N_{\rm tot}^{3/2}\log N_{\rm tot})$ \\
Cost of the first solve:           & $O(N_{\rm tot}^{2})$ \\
Cost of subsequent solves:         & $O(N_{\rm tot}^{3/2})$
\end{tabular}
\end{center}
We observe that even though the asymptotic cost of forming the linear
systems is less than the cost of factoring the matrices, the set-up
phase tends to dominate unless $N_{\rm tot}$ is large. Moreover, the $O(N_{\rm tot}^{3/2})$
cost of the subsequent solves has a very small constant of proportionality.

\lsp

\item
The high order discretization employed achieves high accuracy with a small number
of points. In practical terms, this means that despite the $O(N_{\rm tot}^{2})$
scaling, the scheme is very fast even for moderately complicated geometries.

\lsp

\item The system matrices $\mtx{I} + \mtx{A}^{(n)}$ often have internal structure that allow
them to be inverted using ``fast methods'' such as, e.g., those in
\cite{Martinsson:04a}. The cost of inversion and application can in fact be
accelerated to near optimal complexity.
\end{itemize}


\section{Accelerations for the Single and Double Layer Kernels Associated with Laplace's Equation}
\label{sec:laplace}

This section describes an efficient technique based on recursion relations for
evaluating the kernel $k_{n}$, cf.~(\ref{eq:def_kn}), when $k$ is either the single or
double layer kernel associated with Laplace's equation.

\subsection{The Double Layer Kernels of Laplace's Equation}
\label{sec:doublelayer}

Let $D \subset \mathbb{R}^{3}$ be a bounded domain whose boundary is
given by a smooth surface $\Gamma$, let $E = \bar{D}^{\rm c}$ denote the domain exterior to
$D$, and let $\nb$ be the outward unit normal to $D$. Consider the interior and exterior Dirichlet problems
of potential theory \cite{Guenther:88a},
\begin{align}
    & \Delta u = 0 \hspace{0.5em} \textrm{in} \hspace{0.5em} D, \hspace{1em} u = f \hspace{0.5em} \textrm{on} \hspace{0.5em} \Gamma,
    \hspace{2em} &\textrm{(interior Dirichlet problem)} \label{eq:intDir} \\
    & \Delta u = 0 \hspace{0.5em} \textrm{in} \hspace{0.5em} E, \hspace{1em} u = f \hspace{0.5em} \textrm{on} \hspace{0.5em} \Gamma.
    \hspace{2em} &\textrm{(exterior Dirichlet problem)} \label{eq:extDir}
\end{align} The solutions to (\ref{eq:intDir}) and (\ref{eq:extDir}) can be written
in the respective forms
\begin{align}
\label{eq:laplace_ansatz}
    & u(\xb) = \int_{\Gamma} \frac{\nb(\xb') \cdot (\xb - \xb')}{4 \pi |\xb - \xb'|^{3}} \sigma(\xb') \, dA(\xb'),
    \hspace{1em} \xb \in D,  \\
    & u(\xb) = \int_{\Gamma} \( -\frac{\nb(\xb') \cdot (\xb - \xb')}{4 \pi |\xb - \xb'|^{3}} + \frac{1}{4 \pi |\xb - \xb_{0}|} \) \sigma(\xb') \, dA(\xb'),
    \hspace{1em} \xb \in E,
\end{align}
where $\sigma$ is a boundary charge distribution that can be determined using the boundary
conditions. The point $\xb_{0}$ can be placed at any suitable location in $D$.
The resulting equations are
\begin{align}
    -\frac{1}{2} \sigma(\xb) &+ \int_{\Gamma}
    \frac{\nb(\xb') \cdot (\xb - \xb')}{4 \pi |\xb - \xb'|^{3}} \sigma(\xb') \, dA(\xb') = f(\xb), \label{eq:intDirBie} \\
    -\frac{1}{2} \sigma(\xb) &+ \int_{\Gamma} \( -\frac{\nb(\xb') \cdot (\xb - \xb')}{4 \pi |\xb - \xb'|^{3}} +
    \frac{1}{4 \pi |\xb - \xb_{0}|} \) \sigma(\xb') \, dA(\xb') = f(\xb), \label{eq:extDirBie}
\end{align}
where $\xb \in \Gamma$ in (\ref{eq:intDirBie}) and (\ref{eq:extDirBie}).
\begin{remark}There are other integral formulations for the solution to Laplace's equation.
The double layer formulation presented here is a good choice in that it provides an
integral operator that leads to well conditioned linear systems.  However, the
methodology of this chapter is equally applicable to single-layer formulations that
lead to first kind Fredholm BIEs.\end{remark}


\subsection{Separation of Variables}

Using the procedure given in Section \ref{sec:sepVar}, if $\Gamma = \gamma \times \mathbb{T}$, then (\ref{eq:intDir}) and (\ref{eq:extDir}) can be recast as a series of BIEs defined along $\gamma$.  We express $\nb$ in cylindrical coordinates as
\begin{equation*}
    \nb(\xb') = (n_{r'} \cos \theta', n_{r'} \sin \theta', n_{z'}).
\end{equation*}
Further,
\begin{align*}
    |\xb - \xb'|^{2} &= (r \cos \theta - r' \cos \theta')^{2} + (r \sin \theta - r' \sin \theta')^{2} + (z - z')^2 \\
    &= r^{2} + (r')^{2} - 2 r r' (\sin \theta \sin \theta' + \cos \theta \cos \theta') + (z - z')^2 \\
    &= r^{2} + (r')^{2} - 2 r r' \cos(\theta - \theta') + (z - z')^{2}
\end{align*}
and
\begin{align*}
    \nb(\xb') \cdot (\xb - \xb') &= (n_{r'} \cos \theta', n_{r'} \sin \theta', n_{z'}) \cdot
    (r \cos \theta - r' \cos \theta', r \sin \theta - r' \sin \theta', z - z') \\
    &= n_{r'} r( \sin \theta \sin \theta' + \cos \theta \cos \theta') - n_{r'} r' + n_{z'}(z - z') \\
    &= n_{r'} (r \cos(\theta - \theta') - r') + n_{z'}(z - z').
\end{align*}
Then for a point $\xb' \in \Gamma$, the kernel of the internal Dirichlet problem can be expanded as
\begin{equation*}
    \frac{\nb(\xb') \cdot (\xb - \xb')}{4 \pi |\xb - \xb'|^{3}} = \frac{1}{\sqrt{2 \pi}} \sum_{n \in \mathbb{Z}} e^{i n (\theta - \theta')} d^{(i)}_{n}(r,z,r',z'),
\end{equation*}
where
\begin{equation*}
    d^{(i)}_{n}(r,z,r',z') = \frac{1}{\sqrt{32 \pi^{3}}} \int_{\mathbb{T}} e^{-i n \theta} \left[
     \frac{n_{r'} ( r \cos \theta - r' ) + n_{z'}(z - z')}
        {(r^{2} + (r')^{2} - 2 r r' \cos \theta + (z - z')^{2})^{3/2}} \right] \, d \theta.
\end{equation*}
Similarly, the kernel of the external Dirichlet problem can be written as
\begin{equation*}
    -\frac{\nb(\xb') \cdot (\xb - \xb')}{4 \pi |\xb - \xb'|^{3}} + \frac{1}{4 \pi | \xb - \xb_{0} |} =
    \frac{1}{\sqrt{2 \pi}} \sum_{n \in \mathbb{Z}} e^{i n (\theta - \theta')} d^{(e)}_{n}(r,z,r',z'),
\end{equation*}
with
\begin{align*}
    d^{(e)}_{n}(r,z,r',z') = \frac{1}{\sqrt{32 \pi^{3}}} & \int_{\mathbb{T}} e^{-i n \theta} \biggl(
     -\frac{n_{r'} ( r \cos \theta - r' ) + n_{z'}(z - z')}
        {(r^{2} + (r')^{2} - 2 r r' \cos \theta + (z - z')^{2})^{3/2}} + \\
    & + \frac{1}{(r^{2} + r_{0}^{2} - 2 r r_{0} \cos \theta + (z - z_{0})^{2})^{1/2}} \biggr)
     \, d \theta,
\end{align*}
where $\xb_{0}$ has been written in cylindrical coordinates as $(r_{0} \cos(\theta_{0}),r_{0} \sin(\theta_{0}),z_{0})$.
With the expansions of the kernels available, the procedure described in Section \ref{sec:genAlg} can be used to solve (\ref{eq:intDirBie}) and (\ref{eq:extDirBie}) by solving
\begin{equation}
    \sigma_{n}(r,z) + \sqrt{2 \pi} \int_{\gamma} d^{(i)}_{n}(r,r',z,z') \sigma_{n}(r',z') \, r' \, dl(r',z') = f_{n}(r,z)
    \label{eq:intDirFour}
\end{equation}
and
\begin{equation}
    \sigma_{n}(r,z) + \sqrt{2 \pi} \int_{\gamma} d^{(e)}_{n}(r,r',z,z') \sigma_{n}(r',z') \, r' \, dl(r',z') = f_{n}(r,z),
    \label{eq:extDirFour}
\end{equation}
respectively for $n = -N,-N+1,\ldots,N$.
Note that the kernels $d^{(i)}_{n}$ and $d^{(e)}_{n}$ contain a log-singularity
as $(r',z') \rightarrow (r,z)$.


\subsection{Evaluation of Kernels}
\label{sec:recursion}

The values of $d^{(i)}_{n}$ and $d^{(e)}_{n}$ for $n = -N,-N+1,\ldots,N$ need to be computed efficiently and with high accuracy to construct the Nystr\"{o}m discretization of (\ref{eq:intDirFour}) and (\ref{eq:extDirFour}).  Note that the integrands of $d^{(i)}_{n}$ and $d^{(e)}_{n}$ are real valued and even functions on the interval $[-\pi,\pi]$.  Therefore, $d^{(i)}_{n}$ can be written as
\begin{equation}
    d_{n}^{(i)}(r,z,r',z') = \frac{1}{\sqrt{32 \pi^{3}}} \int_{\mathbb{T}} \left[ \frac{n_{r'} ( r \cos t - r' ) + n_{z'}(z - z')}
        {(r^{2} + (r')^{2} - 2 r r' \cos t + (z - z')^{2})^{3/2}} \right] \cos(n t) \, dt.
    \label{eq:doubleCoef}
\end{equation}
Note that $d^{(e)}_{n}$ can be written in a similar form.

This integrand is oscillatory and increasingly peaked at the origin as
$(r',z')$ approaches $(r,z)$.  As long as $r'$ and $r$ as well as $z'$ and $z$ are
well separated, the integrand does not experience peaks near the origin, and as
mentioned before, the FFT provides a fast and accurate way for calculating
$d^{(i)}_{n}$ and $d^{(e)}_{n}$.

In regimes where the integrand is peaked, the FFT no longer provides a means of
evaluating $d^{(i)}_{n}$ and $d^{(e)}_{n}$ with the desired accuracy.  One possible
solution to this issue is applying adaptive quadrature to fully resolve the peak.
However, this must be done for each value of $n$ required and becomes prohibitively
expensive if $N$ is large.

Fortunately, an analytical solution to (\ref{eq:doubleCoef}) exists.  As noted in
\cite{Cohl:99a}, the single-layer kernel can be expanded with respect to the azimuthal
variable as
\begin{align*}
    s(\xb,\xb') = \frac{1}{4 \pi | \xb - \xb' |} &= \frac{1}{4 \pi (r^{2} + (r')^{2} - 2 r r' \cos(\theta - \theta') + (z - z')^{2})^{1/2}} \\
    &= \frac{1}{\sqrt{2 \pi}} \sum_{n \in \mathbb{Z}} e^{i n (\theta - \theta')} s_{n}(r,z,r',z'),
\end{align*}
where
\begin{align*}
    s_{n}(r,z,r',z') &= \frac{1}{\sqrt{32 \pi^{3}}} \int_{\mathbb{T}}
    \frac{\cos( n t)}{(r^{2} + (r')^{2} - 2 r r' \cos(t) + (z - z')^{2})^{1/2}} \, dt \\
    &= \frac{1}{ \sqrt{8 \pi^{3} r r'}} \int_{\mathbb{T}} \frac{\cos(nt)}{\sqrt{8 (\chi - \cos(t))} } \, dt \\
    &= \frac{1}{\sqrt{8 \pi^{3} r r'}} \mathcal{Q}_{n-1/2}(\chi),
\end{align*}
$\mathcal{Q}_{n-1/2}$ is the half-integer degree Legendre function of the second kind, and
\begin{equation*}
    \chi = \frac{r^{2} + (r')^{2} + (z - z')^{2}}{2 r r'}.
\end{equation*}

To find an analytical form for (\ref{eq:doubleCoef}), first note that in cylindrical coordinates the double-layer kernel can be written in terms of the single-layer kernel,
\begin{align*}
    \frac{\nb (\xb') \cdot (\xb - \xb')}{4 \pi |\xb - \xb'|^{3}} =& \frac{n_{r'} ( r \cos(\theta - \theta') - r' ) + n_{z'}(z - z')}
        {4 \pi (r^{2} + (r')^{2} - 2 r r' \cos(\theta - \theta') + (z - z')^{2})^{3/2}} \\
    =& \frac{1}{4 \pi} \biggl[ n_{r'} \frac{\p}{\p r'} \( \frac{1}{(r^{2} +
    (r')^{2} - 2 r r' \cos(\theta - \theta') + (z - z')^{2})^{1/2}} \) + \\
    &+ n_{z'} \frac{\p}{\p z'} \( \frac{1}{(r^{2} + (r')^{2} - 2 r r' \cos(\theta - \theta') + (z - z')^{2})^{1/2}} \) \biggr].
\end{align*}
The coefficients of the Fourier series expansion of the double-layer kernel are then given by $d^{(i)}_{n}$, which can be written using the previous equation as
\begin{align*}
    d_{n}^{(i)}(r,z,r',z') =& n_{r'} \int_{\mathbb{T}} \frac{\p}{\p r'} \( \frac{\cos(n t)}{ (32 \pi^{3} (r^{2} +
    (r')^{2} - 2 r r' \cos(t) + (z - z')^{2}))^{1/2}} \) \, dt + \\
    &+ n_{z'} \int_{\mathbb{T}} \frac{\p}{\p z'}
    \( \frac{\cos(n t)}{(32 \pi^{3}(r^{2} + (r')^{2} - 2 r r' \cos(t) + (z - z')^{2}))^{1/2}} \)  \, dt \\
    =& n_{r'} \frac{\p}{\p r'} \( \frac{1}{\sqrt{8 \pi^{3} r r'}} \mathcal{Q}_{n-1/2}(\chi) \) +
    n_{z'} \frac{\p}{\p z'} \( \frac{1}{\sqrt{8 \pi^{3} r r'}} \mathcal{Q}_{n-1/2}(\chi) \) \\
    =& \frac{1}{\sqrt{8 \pi^{3} r r'}}\[ n_{r'} \( \frac{\p \mathcal{Q}_{n-1/2}(\chi)}{\p \chi} \frac{\p \chi}{\p r'} -
    \frac{\mathcal{Q}_{n-1/2}(\chi)}{2 r'} \)
    + n_{z'} \frac{\p \mathcal{Q}_{n-1/2}(\chi)}{\p \chi} \frac{\p \chi}{\p z'} \].
\end{align*}
To utilize this form of $d^{(i)}_{n}$, set $\mu = \sqrt{\frac{2}{\chi + 1}}$ and note that
\begin{align*}
    &\frac{\p \chi}{\p r'} = \frac{(r')^{2} - r^{2} - (z - z')^{2}}{2 r (r')^{2}}, \\
    &\frac{\p \chi}{\p z'} = \frac{z' - z}{r r'}, \\
    &\mathcal{Q}_{-1/2}(\chi) = \mu K(\mu), \\
    &\mathcal{Q}_{1/2}(\chi) = \chi \mu K(\mu) - \sqrt{2(\chi + 1)} E(\mu), \\
    &\mathcal{Q}_{-n-1/2}(\chi) = \mathcal{Q}_{n-1/2}(\chi), \\
    &\mathcal{Q}_{n-1/2}(\chi) = 4 \frac{n-1}{2n-1} \chi \mathcal{Q}_{n-3/2}(\chi) - \frac{2n-3}{2n-1}\mathcal{Q}_{n-5/2}(\chi), \\
    &\frac{\p \mathcal{Q}_{n-1/2}(\chi)}{\p \chi} = \frac{2n-1}{2(\chi^2-1)} \( \chi \mathcal{Q}_{n-1/2} - \mathcal{Q}_{n-3/2} \),
\end{align*}
where $K$ and $E$ are the complete elliptic integrals of the first and second kinds, respectively.  The first two relations follow immediately from the definition of $\chi$ and the relations for the Legendre functions of the second kind can be found in \cite{Abramowitz:65a}.  With these relations in hand, the calculation of $d^{(i)}_{n}$ for $n = -N,-N+1,\ldots,N$ can be done accurately and efficiently when $r'$ and $r$ as well as $z'$ and $z$ are in close proximity.  The calculation of $d^{(e)}_{n}$ can be done analogously.

\begin{remark}Note that the forward recursion relation for the Legendre functions $\mathcal{Q}_{n-1/2}(\chi)$ is unstable when $\chi > 1$.  In practice, the instability is mild when $\chi$ is near $1$ and the recursion relation can still be employed to accurately compute values in this regime.  Additionally, if stability becomes an issue, Miller's algorithm \cite{Gil:07a} can be used to calculate the values of the Legendre functions using the backwards recursion relation, which is stable for $\chi > 1$.\end{remark}


\section{Fast Kernel Evaluation for the Helmholtz Equation}
\label{sec:kerExtension}

Section \ref{sec:laplace} describes how to efficiently evaluate the
kernels $k_{n}$ as defined by (\ref{eq:def_kn}) for kernels associated
with Laplace's equation. This section generalizes these methods
to a broad class of kernels that includes the single and double layer
kernels associated with the Helmholtz equation.

\subsection{Rapid Kernel Calculation via Convolution}
\label{sec:kerExt}

Consider a kernel of the form
\begin{equation}
	f(\xb,\xb') = s(\xb,\xb') \, g(\xb,\xb'),
	\label{eq:ch4:prodker}
\end{equation}
where
\begin{align*}
    s(\xb,\xb') = \frac{1}{4 \pi | \xb - \xb' |}
\end{align*}
is the single layer kernel of Laplace's equation and $g(\xb,\xb')$ is a smooth function for all $\xb,\xb' \in \mathbb{R}^{3}$.  Common examples of kernels that take this form include the fundamental solution of the Helmholtz equation and screened Coulomb (Yukawa) potentials.

Letting
\begin{align*}
	\xb  &= (r\,\cos\theta,\,r\,\sin\theta,\,z),\\
	\xb' &= (r'\,\cos\theta',\,r'\,\sin\theta',\,z'),
\end{align*}
we are interested in calculating the Fourier expansion of (\ref{eq:ch4:prodker}) in terms of the azimuthal variable.  When $ g(\xb,\xb') = 1$ (the Laplace kernel), we know how to rapidly compute these Fourier coefficients rapidly and efficiently.  However, when $g$ takes a nontrivial form, this is not generally true; there is no known analytical formula for calculating the Fourier coefficients of (\ref{eq:ch4:prodker}).

We will now describe an efficient technique for calculating the the Fourier coefficients of (\ref{eq:ch4:prodker}), when the function $g$ is sufficiently smooth.  For a fixed value of $(r,z)$ and $(r',z')$, the functions $s$ and $g$ are periodic in the azimuthal variable over the interval $\mathbb{T}$.  Dropping the dependence of $s$ and $g$ on $(r,z)$ and $(r',z')$ for notational clarity, we define $t = \theta-\theta' \in \mathbb{T}$ and the Fourier series expansions of $s$ and $g$ as
\begin{align}
	s(t) &= \sum_{n \in \mathbb{Z}}\frac{e^{int}}{\sqrt{2\pi}}\,s_{n}, \label{eq:sers} \\
	g(t) &= \sum_{n \in \mathbb{Z}}\frac{e^{int}}{\sqrt{2\pi}}\,g_{n}, \label{eq:serg}
\end{align}
where
\begin{align}
    s_{n} &= \int_{\mathbb{T}} \frac{e^{-int}}{\sqrt{2\pi}}\, s(t) \, dt,\label{eq:exps} \\
    g_{n} &= \int_{\mathbb{T}} \frac{e^{-int}}{\sqrt{2\pi}}\, g(t) \, dt. \label{eq:expg}
\end{align}
The values given by (\ref{eq:exps}) can be calculated as described in Section \ref{sec:laplace}, while the values given by (\ref{eq:expg}) can be rapidly and accurately computed using the FFT.

Assuming the Fourier series defined by (\ref{eq:sers}) and (\ref{eq:serg}) are uniformly convergent, we find
$$
	f_{n} = \frac{1}{\sqrt{2 \pi}} \int_{\mathbb{T}} s(t) \, g(t) \, e^{- i n t} \, dt
	= \sum_{k \in \mathbb{Z}} s_{k} \, g_{n - k}
	= [s_{k} * g_{k}](n),
$$
where $s_{k} * g_{k}$ is the discrete convolution of the sequences defined by (\ref{eq:exps}) and (\ref{eq:expg}).

In a practical setting, the Fourier series are truncated to finite length.  Assuming that
we have kept $-N,-N+1,\ldots,N$ terms, directly calculating the convolution
would require $O(N^{2})$ operations.   Fortunately, this computation can be accelerated
to $O(N \log N)$ operations by employing the discrete convolution theorem and the
FFT \cite{Briggs:95a}.

Letting $\mathcal{D}$ denote the discrete Fourier transform (DFT), the discrete convolution theorem states that the convolution of two periodic sequences $\{a_{n}\}$ and $\{b_{n}\}$ is related by
\begin{equation*}
	\mathcal{D}\{ a_{n} * b_{n} \}_{k} = \alpha A_{k} B_{k},
\end{equation*}
where $\{A_{n}\} = \mathcal{D} \{a_{n}\}$, $\{B_{n}\} = \mathcal{D} \{b_{n}\}$, and $\alpha$ is a known constant depending upon the length of the periodic sequences and the definition of the DFT.  Thus, we can rapidly calculate the convolution of two periodic sequences by taking the FFT of each sequence, computing the pointwise product of the result, and then applying the inverse FFT.

Of course, the sequences that we need to convolute are not periodic.  Applying the discrete convolution to the sequences defined by (\ref{eq:exps}) and (\ref{eq:expg}) will not be exact, but the error incurred will be small assuming that the Fourier coefficients decay rapidly and that $N$ is large enough.  To see this, assume that
\begin{align*}
	s(t) &= \sum_{n = -N}^{N} \frac{e^{int}}{\sqrt{2\pi}}\,s_{n}, \\
	g(t) &= \sum_{n = -N}^{N} \frac{e^{int}}{\sqrt{2\pi}}\,g_{n}.
\end{align*}
Then the exact Fourier representation of $f$ can be found by taking the
product of these two series, which will be of length $4N+1$.  As is well known,
the coefficients of this product is given by the discrete convolution of the
sequences containing the coefficients of the two series, and these sequences must
first be padded with $2 N$ zeros.  Thus, we can effectively calculate the
Fourier coefficients of the function given by (\ref{eq:ch4:prodker}) by calculating
$2N+1$ Fourier coefficients of $s$ and $g$, padding these sequences with zeros,
calculating the discrete convolution of these sequences, and truncating the
resulting sequence.  In practice, padding may not even be required if the Fourier
coefficients of $s$ and $g$ decay sufficiently fast.

Note that the procedure described in this section is quite general.  The azimuthal
Fourier coefficients of many kernels that can be represented as the product of a
singular function and a smooth function can be found, assuming that there is an
accurate technique for determining the coefficients of the singular function.

\subsection{Application to the Helmholtz Equation}
\label{sec:helmholtz}

In this section, we will apply the fast kernel calculation technique described in Section \ref{sec:kerExt} to the exterior Dirichlet problem for the Helmholtz equation.  Let $D \subset \mathbb{R}^{3}$ be a bounded domain whose boundary is given by a smooth surface $\Gamma$, let $E = \bar{D}^{\rm c}$ denote the domain exterior to $D$, and let $\nb$ and be the outward unit normal to $D$.  The partial differential equation representing this problem is given by
\begin{equation}
	\Delta u + k^{2} u = 0 \hspace{0.5em} \textrm{in} \hspace{0.5em} E, \hspace{1em} u =
	f \hspace{0.5em} \textrm{on} \hspace{0.5em} \Gamma, \label{eq:extDirHelm}
\end{equation}
where $k > 0$ is the wavenumber, and $u$ satisfies the Sommerfeld radiation condition
\begin{equation}
	\lim_{r \rightarrow \infty} r \( \frac{\p u}{\p r} - i \, k \, u \) = 0, \label{eq:sommerfeld}
\end{equation}
where $r = |\xb|$ and the limit holds uniformly in all directions $\xb / |\xb|$. Let the single and double layer potentials for the Helmholtz equation be given by
\begin{align}
	\phi(\xb,\xb') &= \frac{e^{i k |\xb - \xb'|}}{4 \pi |\xb - \xb'|},  \hspace{2em} (\textrm{single layer}) \label{eq:singleLayerHelm} \\
	\frac{\p \phi(\xb,\xb')}{\p \nb(\xb')}  &= \frac{\nb(\xb') \cdot (\xb - \xb')}{4 \pi |\xb - \xb'|^{3}} \[ \(1 - i k |\xb - \xb'| \) e^{i k |\xb - \xb'|} \].
	\hspace{2em} (\textrm{double layer})
\label{eq:doubleLayerHelm}
\end{align}
The solution to (\ref{eq:extDirHelm}) can be written in terms of the double layer potential,
\begin{equation*}
	u(\xb) = \int_{\Gamma} \frac{\p \phi(\xb,\xb')}{\p \nb(\xb')} \, \sigma(\xb') \, dA(\xb'),
	\hspace{1em} \xb \in E,
\end{equation*}
where $\sigma$ is a boundary charge density that can be determined using the boundary conditions. The resulting boundary integral equation is given by
\begin{equation}
	\frac{1}{2} \sigma(\xb) + \int_{\Gamma} \frac{\p \phi(\xb,\xb')}{\p \nb(\xb')} \, \sigma(\xb') \, dA(\xb') = f(\xb).
	\label{eq:extDirBieHelm} \end{equation} As is well known,
(\ref{eq:extDirBieHelm}) is not always uniquely solvable, even though
(\ref{eq:extDirHelm}) is uniquely solvable for all $k > 0$.  A common solution to
this is to represent  the solution to (\ref{eq:extDirHelm}) as a combined single and
double layer potential,
\begin{equation*}
	u(\xb) = \int_{\Gamma} \( \frac{\p \phi(\xb,\xb')}{\p \nb(\xb')} - i \, \nu \, \phi(\xb,\xb') \) \, \sigma(\xb') \, dA(\xb'),
	\hspace{1em} \xb \in E, \end{equation*}
where $\nu > 0$.  We have freedom in
choosing $\nu$, see, e.g., \cite{Bruno:01a,Kress:83a}, for some analysis on
this choice.  The boundary integral equation we need to solve is
\begin{equation}
	\frac{1}{2} \sigma(\xb) + \int_{\Gamma} \( \frac{\p \phi(\xb,\xb')}{\p \nb(\xb')} - i \, \nu \, \phi(\xb,\xb') \) \, \sigma(\xb') \, dA(\xb') = f(\xb).
	\label{eq:extDirBieHelmComb}
\end{equation}


\section{Fast evaluation of fundamental solutions in cylindrical coordinates}
\label{sec:fundamentalsolution}

The techniques for kernel evaluations described in
Sections \ref{sec:fastform_An}, \ref{sec:laplace}, and \ref{sec:kerExtension} are
useful not only for solving BIEs, but for solving the Laplace
and Helmholtz equations in a variety of contexts where
cylindrical coordinates are effective. To illustrate, observe
that the free space equation
\begin{equation}
\label{eq:freespace}
-\Delta u(\xb) - k^{2}\,u(\xb) = f(\xb),\qquad\xb \in \mathbb{R}^{3}
\end{equation}
has the solution
\begin{equation}
\label{eq:freespacesoln}
u(\xb) = \int_{\mathbb{R}^{3}}\phi^{(k)}(\xb,\xb')\,f(\xb')\,dA(\xb').
\qquad \xb \in \mathbb{R}^{3},
\end{equation}
where $\phi^{(k)}$ is the fundamental solution
$$
\phi^{(k)}(\xb) = \frac{e^{ik|\xb|}}{4\pi|\xb|}.
$$
In cylindrical coordinates, we write (\ref{eq:freespacesoln}) as
\begin{equation}
\label{eq:freespacesoln_cyl}
u(\xb) =
\sum_{n=-\infty}^{\infty}
\frac{e^{in\theta}}{\sqrt{2\pi}}
\int_{H}\phi_{n}^{(k)}(r,z,r',z')f_{n}(r',z')dA(r',z')
\end{equation}
where $H = \{(r,z) \in \mathbb{R}^{2}\,\colon\,r \geq 0\}$ is a half-plane,
and where $u_{n}$, $f_{n}$, and $\phi_{n}^{(k)}$ are the Fourier coefficients
defined by
\begin{align*}
u(\xb) =& \sum_{n=-\infty}^{\infty} \frac{e^{in\theta}}{\sqrt{2\pi}} u_{n}(r,z),\\
f(\xb) =& \sum_{n=-\infty}^{\infty} \frac{e^{in\theta}}{\sqrt{2\pi}} f_{n}(r,z),\\
\phi^{(k)}(\xb,\xb') &= \sum_{n=-\infty}^{\infty} \frac{e^{in(\theta-\theta')}}{\sqrt{2\pi}}
\phi_{n}^{(k)}(r,z,r',z').
\end{align*}
The kernel $\phi_{n}^{(k)}$ in (\ref{eq:freespacesoln_cyl})
can be evaluated efficiently using the techniques of Sections \ref{sec:fastform_An}, \ref{sec:laplace}
and \ref{sec:kerExtension}. If $f$ has a rapidly convergent Fourier series, and
if ``fast'' summation (e.g.~the Fast Multipole Method) is used to evaluate the
integrals in (\ref{eq:freespacesoln_cyl}), then very efficient solvers result.

More generally, we observe that the equation (\ref{eq:freespace}) can be expressed
\begin{equation}
\label{eq:freespace_cyl}
- \frac{\p^{2}u_{n}}{\p^{2}r}
- \frac{1}{r}\frac{\p u_{n}}{\p r}
- \frac{\p^{2}u_{n}}{\p^{2}z}
+ \left(\frac{n^{2}}{r^{2}} - k^{2}\right)u_{n} = f_{n},
\qquad n \in \mathbb{Z}.
\end{equation}
and that the function $\phi_{n}^{(k)}$ is the Green's function of (\ref{eq:freespace_cyl}).


\section{Numerical Results}
\label{sec:num}

\begin{figure}[!t] \begin{center} \begin{minipage}{0.31\linewidth} \begin{flushleft}
\hspace{10mm} \includegraphics[height=.85\linewidth]{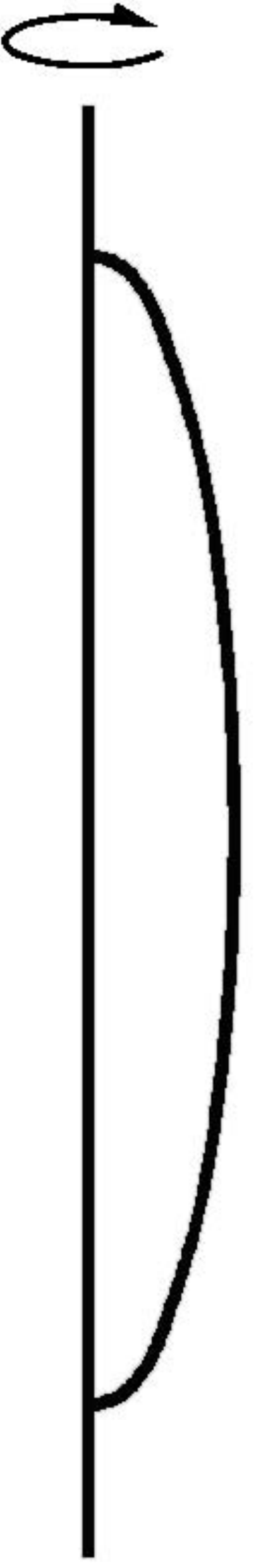}
\end{flushleft} \end{minipage} \begin{minipage}{0.31\linewidth} \begin{center}
\includegraphics[height =.85\linewidth]{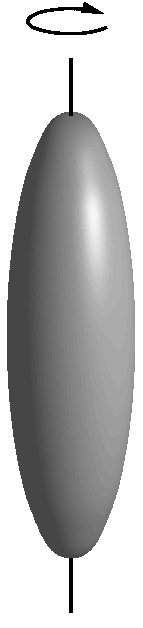} \end{center}
\end{minipage} \\ \begin{minipage}{1\linewidth}\begin{center} (a) \end{center}
\end{minipage} \\ \vspace{4mm} \begin{minipage}{0.31\linewidth} \begin{flushleft}
\hspace{10mm} \includegraphics[height=.75\linewidth]{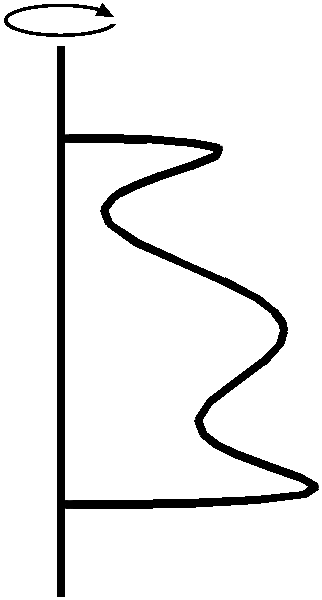}
\end{flushleft} \end{minipage} \begin{minipage}{0.31\linewidth} \begin{center}
\includegraphics[height =.65\linewidth]{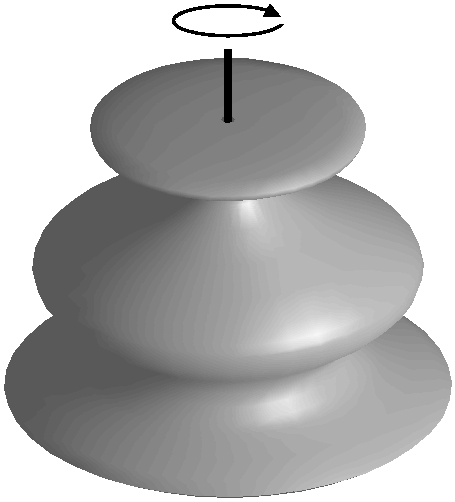} \end{center}
\end{minipage} \\ \begin{minipage}{1\linewidth}\begin{center} (b) \end{center}
\end{minipage} \\ \vspace{4mm} \begin{minipage}{0.31\linewidth} \begin{flushleft}
\hspace{10mm} \includegraphics[height=.75\linewidth]{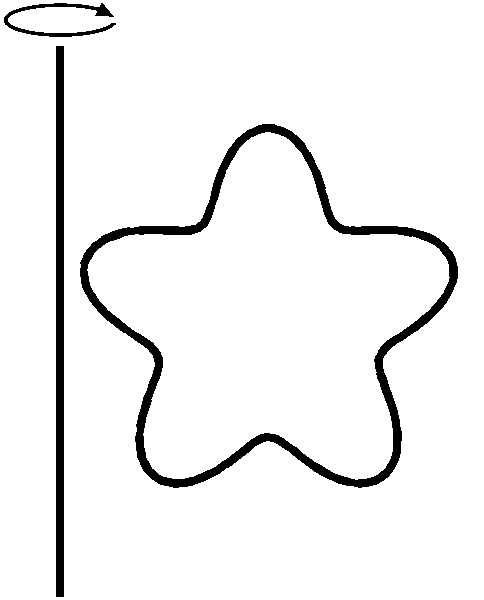}
\end{flushleft} \end{minipage} \begin{minipage}{0.31\linewidth} \begin{center}
\includegraphics[height =.55\linewidth]{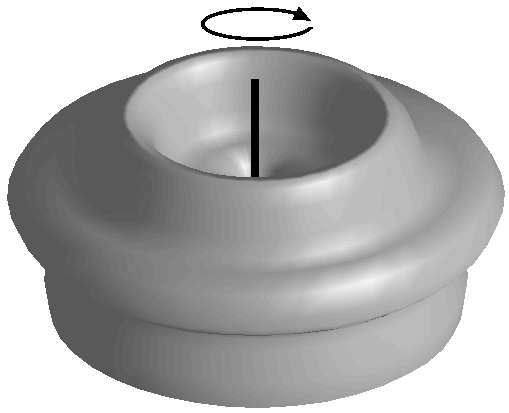} \end{center}
\end{minipage} \\ \begin{minipage}{1\linewidth}\begin{center} (c) \end{center}
\end{minipage} \\ \caption{Domains used in numerical examples. All items are rotated
about the vertical axis.  (a) An ellipse.  (b) A wavy block.  (c) A starfish torus.
} \label{fig:domains} \end{center} \end{figure}

This section describes several numerical experiments performed to assess the
efficiency and accuracy of the numerical scheme outlined in Section
\ref{sec:summary}. The geometries investigated are described in Figure
\ref{fig:domains}. The generating curves were parameterized by arc length, and split
into $N_{\rm P}$ panels of equal length. A 10-point Gaussian quadrature has been
used along each panel, with the modified quadratures of \cite{2001_rokhlin_kolm}
used to handle the integrable
singularities in the kernel.  The algorithm was implemented in FORTRAN, using BLAS,
LAPACK, and the FFT library provided by Intel's MKL library.  All numerical
experiments in this section have been carried out on a Macbook Pro with a 2.4 GHz
Intel Core 2 Duo and 4GB of RAM.

\subsection{Laplace's equation}
\label{sec:numLaplace}

We solved Laplace's equation in the domain interior to the surfaces shown in Figure
\ref{fig:domains}. The solution was represented via the double layer Ansatz
(\ref{eq:laplace_ansatz}) leading to the BIE (\ref{eq:intDirBie}). The kernels
$d_{n}$ defined via (\ref{eq:doubleCoef}) were evaluated using the techniques
described in Section \ref{sec:recursion}. Tn this case $\mtx{A}^{(n)} =
\mtx{A}^{(-n)}$, and so we need only to invert $N + 1$ matrices. Further,
the FFT used here is complex-valued, and a real-valued FFT would yield a significant
decrease in computation time.

To investigate the speed of the proposed method, we solved a sequence of problems on
the domain in Figure \ref{fig:domains}(a). The timing results are given in Table
\ref{tbl:tb1}. The reported results include:

\lsp

\begin{tabular}{l l}
$N_{\rm{P}}$ & the number  of panels used to discretize the contour (each panel has $I/N_{\rm{P}}$ nodes) \\
$N$ & the Fourier truncation parameter (we keep $2N+1$ modes) \\
$T_{\textrm{mat}}$ & time to construct the linear systems (utilizing the recursion relation) \\
$T_{\textrm{inv}}$ & time to invert the linear systems \\
$T_{\textrm{fft}}$ & time to Fourier transform the right hand side and the solution \\
$T_{\textrm{apply}}$ & time to apply the inverse to the right hand side
\end{tabular}

\lsp

The most expensive component of the calculation is the kernel evaluation required to
form the coefficient matrices $\mtx{A}^{(n)}$. Table \ref{tbl:tb2} compares the use
of the recursion relation in evaluating the kernel when it is near-singular to using
an adaptive Gaussian quadrature. The efficiency of the recursion relation is
evident.

\begin{figure}[!t]
\begin{center}
\includegraphics[width =.65\linewidth]{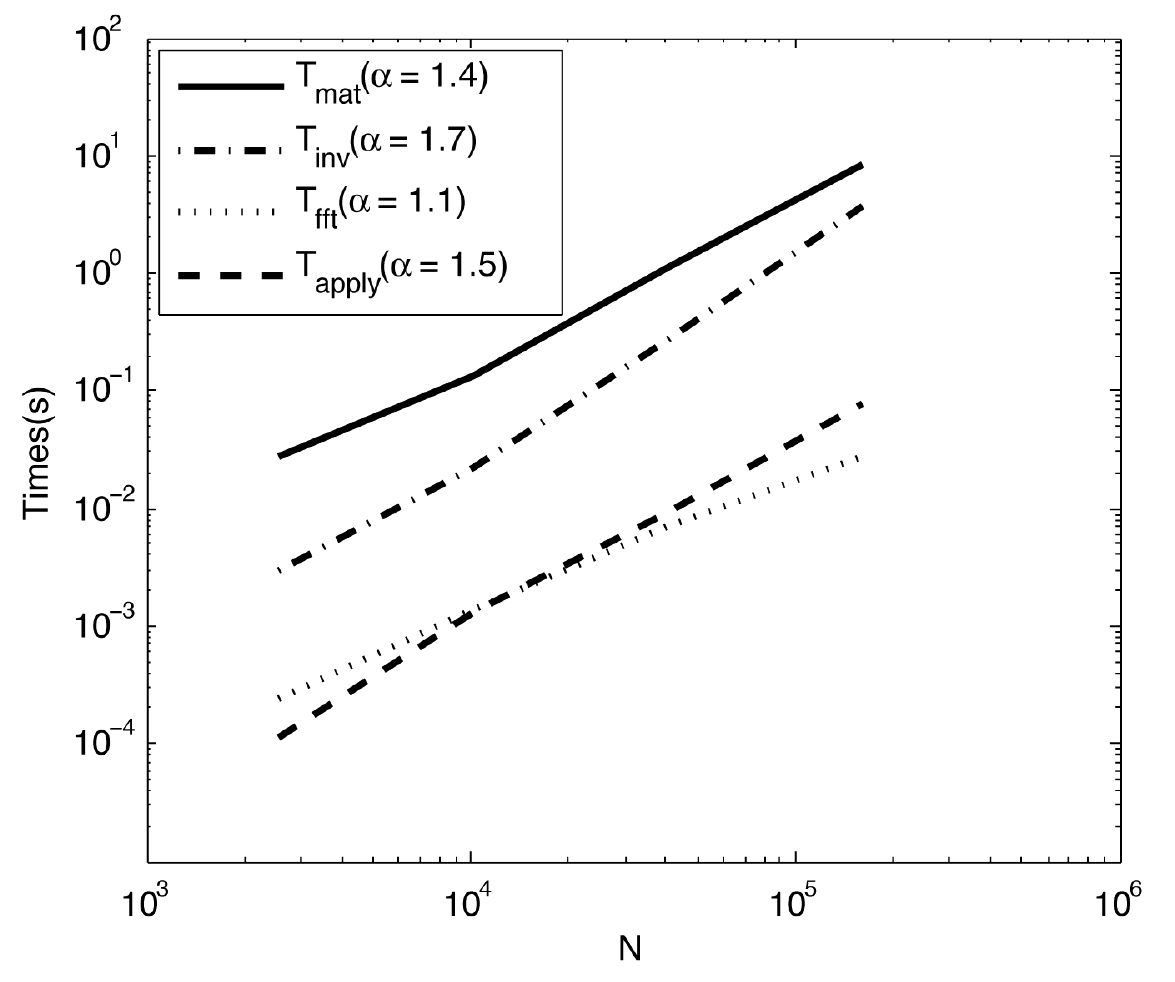}
\end{center}
\caption{Timings of the algorithm as the number of degrees of freedom $N_{\rm tot} = I(2N+1)$ increases.
The timings reported here are for the case $I \approx 2N+1$.
The numbers in parentheses provide estimates of the asymptotic complexity,
i.e.~the best fit to a curve $T = C\,N_{\rm tot}^{\alpha}$. }
\label{fig:scalings}
\end{figure}

Figure \ref{fig:scalings} plots the time to construct the linear systems as the
number of degrees of freedom $N_{\rm tot} = I (2N+1)$ increases, for the case
when $I \approx 2N+1$.  The estimated asymptotic costs given in this
Figure match well with the estimates derived in Section \ref{sec:compcosts}. It is
also clear that as $N_{\rm tot}$ grows, the cost of inversion will eventually dominate. We
remark that the asymptotic scaling of this cost can be lowered by using fast
techniques for the inversion of boundary integral operators, but that little gain
would be achieved for the problem sizes considered here.

We observe that the largest problem reported in Table \ref{tbl:tb1} involves
$320\,800$ degrees of freedom. The method requires $1$ minute of pre-computation for
this example, and is then capable of computing a solution $u$ from a given data
function $f$ in $0.49$ seconds.

\begin{table}[!t]
\begin{center}
\begin{tabular}{| c | c | c | c | c | c | }
\hline
$N_{\rm{P}}$ & $2N + 1$ & $T_{\textrm{mat}}$ & $T_{\textrm{inv}}$ & $T_\textrm{{fft}}$ & $T_{\textrm{apply}}$ \\ \hline
5 & 25 & 1.70E-02 & 1.42E-03 & 7.81E-05 & 3.83E-05 \\
10 & 25 & 3.64E-02 & 6.15E-03 & 1.68E-04 & 2.66E-04 \\
20 & 25 & 9.73E-02 & 3.52E-02 & 3.69E-04 & 2.10E-03 \\
40 & 25 & 3.09E-01 & 2.35E-01 & 6.69E-04 & 4.82E-03 \\
80 & 25 & 1.20E+00 & 1.88E+00 & 1.36E-03 & 2.85E-02 \\ \hline

5 & 51 & 2.83E-02 & 3.02E-03 & 2.38E-04 & 1.13E-04 \\
10 & 51 & 6.71E-02 & 1.23E-02 & 4.48E-04 & 6.58E-04 \\
20 & 51 & 2.02E-01 & 7.42E-02 & 9.17E-04 & 2.63E-03 \\
40 & 51 & 7.28E-01 & 4.94E-01 & 1.92E-03 & 1.03E-02 \\
80 & 51 & 3.07E+00 & 3.73E+00 & 3.59E-03 & 6.17E-02 \\ \hline

5 & 101 & 5.08E-02 & 5.48E-03 & 7.27E-04 & 2.38E-04 \\
10 & 101 & 1.35E-01 & 2.27E-02 & 1.41E-03 & 1.33E-03 \\
20 & 101 & 4.39E-01 & 1.32E-01 & 2.73E-03 & 4.60E-03 \\
40 & 101 & 1.98E+00 & 1.04E+00 & 6.20E-03 & 2.14E-02 \\
80 & 101 & 7.13E+00 & 7.04E+00 & 1.12E-02 & 1.12E-01 \\ \hline

5 & 201 & 1.07E-01 & 1.06E-02 & 1.80E-03 & 6.36E-04 \\
10 & 201 & 3.33E-01 & 4.96E-02 & 3.83E-03 & 2.79E-03 \\
20 & 201 & 1.12E+00 & 2.73E-01 & 7.05E-03 & 9.46E-03 \\
40 & 201 & 4.63E+00 & 1.88E+00 & 1.46E-02 & 4.15E-02 \\
80 & 201 & 1.71E+01 & 1.41E+01 & 2.93E-02 & 2.15E-01 \\ \hline

5 & 401 & 1.87E-01 & 2.15E-02 & 3.43E-03 & 1.48E-03 \\
10 & 401 & 5.85E-01 & 9.51E-02 & 6.59E-03 & 5.05E-03 \\
20 & 401 & 2.15E+00 & 5.42E-01 & 1.33E-02 & 1.91E-02 \\
40 & 401 & 8.40E+00 & 3.71E+00 & 2.78E-02 & 7.98E-02 \\
80 & 401 & 3.15E+01 & 2.83E+01 & 5.56E-02 & 4.34E-01 \\ \hline
\end{tabular}
\vspace{2mm}
\caption{Timing results in seconds performed for the domain given in Figure \ref{fig:domains}(a) for the interior Dirichlet problem.}
\label{tbl:tb1}
\end{center}
\end{table}


\begin{table}[!ht]
\begin{center}
\begin{tabular}{| c | c | c |}
\hline
$2N + 1$ & Composite Quadrature & Recursion Relation \\
\hline
25 & 1.9 & 0.017 \\
50 & 3.1 & 0.028 \\
100 & 6.6 & 0.051 \\
200 & 18.9 & 0.107 \\
\hline
\end{tabular}
\vspace{2mm}
\caption{Timing comparison in seconds for constructing the matrices $(\mtx{I} + \mtx{A}^{(n)})$ using composite Gaussian quadrature and the recursion relation described in Section \ref{sec:recursion}  to evaluate $k_{n}$ for diagonal and near diagonal blocks.  The FFT is used to evaluate $k_{n}$ at all other entries.  $2N + 1$ is the total number of Fourier modes used. 5 panels were used to discretize the boundary.}
\label{tbl:tb2}
\end{center}
\end{table}

To test the accuracy of the approach, we have solved a both interior and exterior
Dirichlet problems on each of the domains given in Figure \ref{fig:domains}. Exact
solutions were generated by placing a few random point charges outside of the domain
where the solution was calculated. The solution was evaluated at points defined on a
sphere encompassing (or interior to) the boundary. The errors reported in Tables
\ref{tbl:tb3}-\ref{tbl:tb5} are relative errors measured in the $l^{\infty}$-norm,
$||u_{\epsilon} - u||_{\infty} / || u ||_{\infty}$, where $u$ is the exact potential
and $u_{\epsilon}$ is the potential obtained from the numerical solution.

For all geometries, 10 digits of accuracy has been obtained from a discretization
involving a relatively small number of degrees of freedom, due to the rapid
convergence of the Gaussian quadrature. This is especially advantageous, as the most
expensive component of the algorithm is the construction of the linear systems, the
majority of the cost being directly related to the number of panels used.  Further,
the number of Fourier modes required to obtain 10 digits of accuracy is on the order
of 100 modes. Although not investigated here, the discretization technique naturally
lends itself to nonuniform refinement in the $rz$-plane, allowing one to resolve
features of the generating curve that require finer resolution.

The number of correct digits obtained as the number of panels and number of Fourier
modes increases eventually stalls.  This is a result of a loss of precision in
determining the kernels, as well as cancelation errors incurred when evaluating
interactions between nearby points.  This is especially prominent with the use of
Gaussian quadratures, as points cluster near the ends of the panels.  If more digits
are required, high precision arithmetic can be employed in the setup phase of the
algorithm.


\begin{table}[!ht]
\begin{center}
\begin{minipage}{0.99\linewidth} \begin{center}
 \begin{tabular}{| c | c | c | c | c | c |}
 \hline
 \multicolumn{1}{|c|}{$N_{\rm P}$} & \multicolumn{5}{c|}{$2N+1$} \\
 \hline
- &25&51&101&201&401    \\
   \hline
5&3.9506E-04&4.6172E-04&4.6199E-04&4.6203E-04&4.6204E-04\\
10&1.3140E-05&1.1091E-08&4.8475E-09&4.8480E-09&4.8481E-09\\
20&1.7232E-05&7.7964E-09&4.7197E-12&4.7232E-12&4.7237E-12\\
40&2.7527E-05&2.7147E-08&2.8818E-14&5.7173E-14&5.7658E-14\\
80&2.118E-05&9.4821E-09&2.1529E-13&2.0392E-13&2.0356E-13\\
  \hline
\end{tabular}
\end{center} \end{minipage} \\
\vspace{2mm}
\caption{Error in internal Dirichlet problem solved on domain (a) in Figure \ref{fig:domains}.}
\label{tbl:tb3}
\end{center}
\end{table}


\begin{table}[!ht]
\begin{center}
\begin{minipage}{0.99\linewidth} \begin{center}
 \begin{tabular}{| c | c | c | c | c | c |}
 \hline
 \multicolumn{1}{|c|}{$N_{\rm P}$} & \multicolumn{5}{c|}{$2N+1$} \\
 \hline
- &25&51&101&201&401    \\
   \hline
5  &8.6992E-04&1.3615E-03&1.3620E-03&1.3621E-03&1.3621E-03 \\
10&2.2610E-04&9.6399E-05&9.6751E-05&9.6751E-05&9.6752E-05\\
20&2.6291E-04&4.6053E-07&2.4794E-07&2.4794E-07&2.4794E-07\\
40&3.1714E-04&2.8922E-07&2.2875E-11&2.3601E-11&2.3605E-11\\
80&3.0404E-04&3.5955E-07&3.3708E-11&3.3138E-11&3.3150E-11\\
  \hline
\end{tabular}
\end{center} \end{minipage} \\
\vspace{2mm}
\caption{Error in external Dirichlet problem solved on domain (b) in Figure \ref{fig:domains}.}
\label{tbl:tb4}
\end{center}
\end{table}


\begin{table}[!ht]
\begin{center}
\begin{minipage}{0.99\linewidth} \begin{center}
 \begin{tabular}{| c | c | c | c | c | c |}
 \hline
 \multicolumn{1}{|c|}{$N_{\rm P}$} & \multicolumn{5}{c|}{$2N+1$} \\
 \hline
- &25&51&101&201&401    \\
   \hline
5  &4.3633E-04 & 7.9169E-05 & 7.8970E-05 & 7.8970E-05 & 7.8971E-05 \\
10&3.9007E-04 & 6.8504E-07 & 2.0274E-08 & 2.0272E-08 & 2.0272E-08\\
20&3.8803E-04 & 6.4014E-07 & 3.2138E-11 & 3.1624E-11 & 3.1625E-11\\
40&3.8456E-04 & 6.4098E-07 & 6.5742E-12 & 3.4529E-12 & 3.4530E-12\\
80&3.9828E-04 & 6.4486E-07 & 6.8987E-12 & 3.1914E-12 & 3.1913E-12\\
  \hline
\end{tabular}
\end{center} \end{minipage} \\
\vspace{2mm}
\caption{Error in external Dirichlet problem solved on domain (c) in Figure \ref{fig:domains}.}
\label{tbl:tb5}
\end{center}
\end{table}

\begin{figure}[!t]
   \centering
   \includegraphics[width=0.65\linewidth]{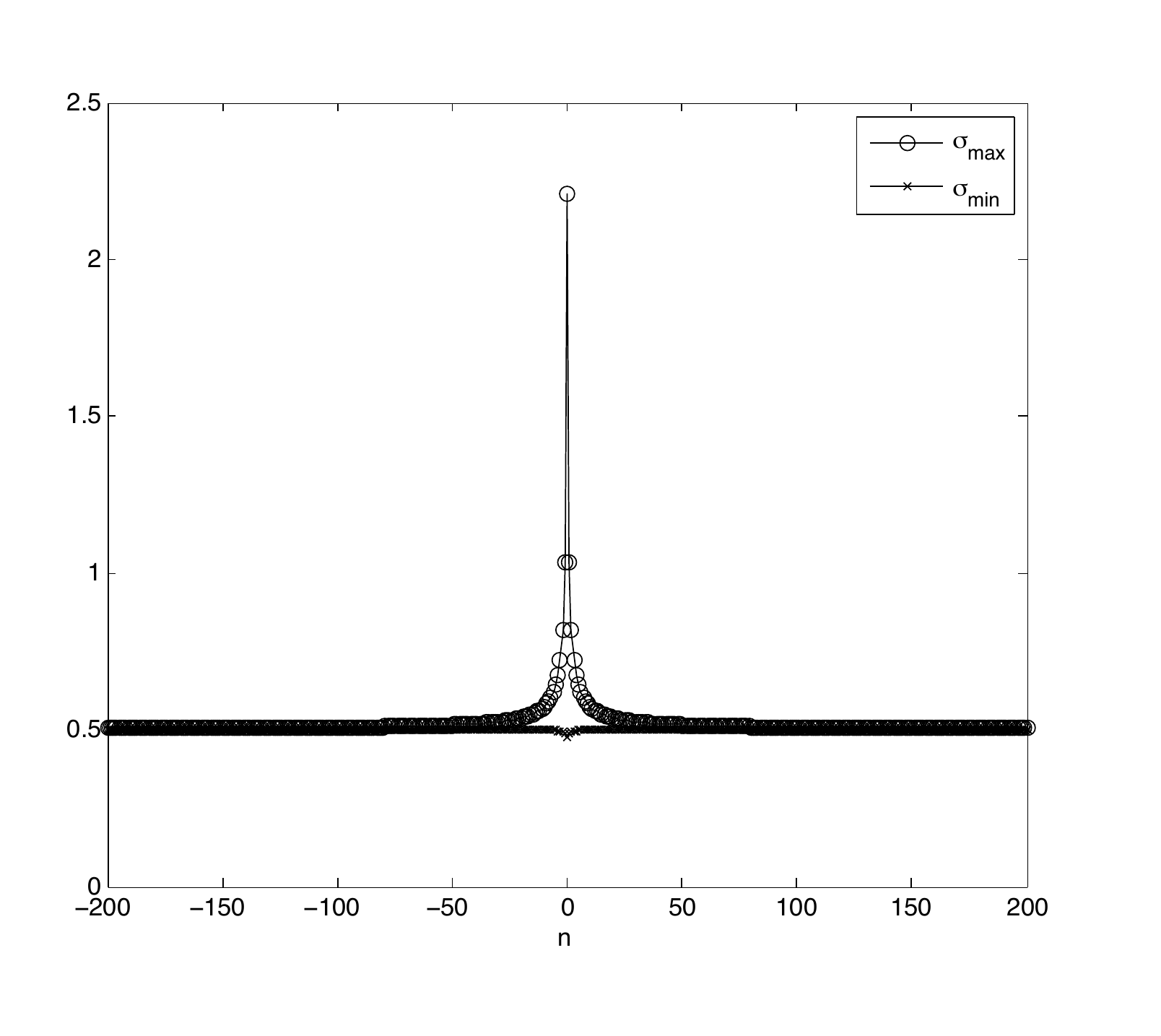}
\caption{Maximum and minimum singular values for the matrices resulting from an 80 panel discretization of a sphere using 400 Fourier modes, where $n$ is the the matrix associated with the $n^{\textrm{th}}$ Fourier mode.}
\label{fig:cond}
\end{figure}

Finally, we investigated the conditioning of the numerical procedure.
Figure \ref{fig:cond} shows the smallest and largest singular values of the matrices
$\{\tfrac{1}{2}\mtx{I} + \mtx{A}^{(n)}\}_{n=-200}^{200}$ on the domain shown in
Figure \ref{fig:domains}(a). The convergence of  both the smallest and the largest
singular values to $1/2$ follow from the convergence $||\mtx{A}^{(n)}|| \rightarrow 0$,
cf.~Section \ref{sec:trunction}. Figure \ref{fig:cond} indicates both
that all matrices involved are well-conditioned, and that truncation of the Fourier
series is generally safe.

\subsection{Helmholtz Equation}
\label{numHelm}

In this section, we repeat many of the experiments reported in Section
\ref{sec:numLaplace}, but now for the Helmholtz equation on an exterior domain with
the associated ``combined field'' BIE formulation (\ref{eq:extDirBieHelmComb}). The
algorithm employed to solve the integral equation is the same as described in
Section \ref{sec:genAlg}, with the caveat that the kernels are calculated using the
fast procedure described in Section \ref{sec:kerExt}.

Table \ref{tbl:helmtime} presents timing results, with all variables defined as in
Section \ref{sec:numLaplace}.

\begin{table}[!ht]
\begin{center}
\begin{tabular}{| c | c | c | c | c | c | }
\hline
$N_{\rm{P}}$ & $2N+1$ & $T_{\textrm{mat}}$ &
$T_{\textrm{inv}}$ & $T_\textrm{{fft}}$ & $T_{\textrm{apply}}$ \\
\hline
5 & 25 & 4.51E-02 & 3.78E-03 & 2.59E-04 & 1.77E-04 \\
10 & 25 & 1.18E-01 & 2.03E-02 & 4.42E-04 & 8.38E-04 \\
20 & 25 & 3.77E-01 & 1.48E-01 & 8.47E-04 & 3.05E-03 \\
40 & 25 & 1.27E+00 & 9.71E-01 & 1.71E-03 & 1.19E-02 \\
80 & 25 & 5.37E+00 & 8.20E+00 & 3.68E-03 & 8.17E-02 \\ \hline

5 & 51 & 8.06E-02 & 6.91E-03 & 5.09E-04 & 3.66E-04 \\
10 & 51 & 2.25E-01 & 4.27E-02 & 1.10E-03 & 1.00E-03 \\
20 & 51 & 7.53E-01 & 2.85E-01 & 2.08E-03 & 6.13E-03 \\
40 & 51 & 2.69E+00 & 1.95E+00 & 3.83E-03 & 2.79E-02 \\
80 & 51 & 1.01E+01 & 1.47E+01 & 7.55E-03 & 1.40E-01 \\ \hline

5 & 101 & 1.57E-01 & 1.36E-02 & 1.33E-03 & 8.13E-04 \\
10 & 101 & 4.74E-01 & 8.20E-02 & 2.55E-03 & 3.05E-03 \\
20 & 101 & 1.58E+00 & 5.27E-01 & 5.03E-03 & 1.17E-02 \\
40 & 101 & 5.95E+00 & 3.81E+00 & 1.01E-02 & 5.67E-02 \\
80 & 101 & 2.11E+01 & 2.72E+01 & 1.90E-02 & 2.39E-01 \\ \hline

5 & 201 & 3.02E-01 & 2.50E-02 & 2.56E-03 & 1.64E-03 \\
10 & 201 & 9.40E-01 & 1.56E-01 & 5.09E-03 & 5.71E-03 \\
20 & 201 & 3.23E+00 & 1.02E+00 & 1.01E-02 & 2.15E-02 \\
40 & 201 & 1.19E+01 & 7.32E+00 & 2.05E-02 & 9.94E-02 \\
80 & 201 & 4.35E+01 & 5.38E+01 & 4.04E-02 & 4.67E-01 \\ \hline

5 & 401 & 5.90E-01 & 5.046E-02 & 5.66E-03 & 3.48E-03 \\
10 & 401 & 1.86E+00 & 2.97E-01 & 1.05E-02 & 1.06E-02 \\
20 & 401 & 6.60E+00 & 2.04E+00 & 2.20E-02 & 4.75E-02 \\
40 & 401 & 2.40E+01 & 1.46E+01 & 4.42E-02 & 1.98E-01 \\ \hline
\end{tabular}
\vspace{2mm}
\caption{Timing results in seconds performed for a spherical domain.}
\label{tbl:helmtime}
\end{center}
\end{table}

The largest problem size considered here has 80 panels and 201 Fourier modes,
leading to $160 \, 800$ unknowns discretizing the surface. Note that this problem
size is slightly smaller than the largest considered in Section
\ref{sec:numLaplace}, due to memory constraints.  This is because the matrices now
contain complex entries, and thus use twice the memory compared with the Laplace
case.  The total running time of the algorithm is 97 seconds for this problem size.
If given additional right hand sides, we can solve them in 0.51 seconds. We also
remark that the asymptotic scaling of the cost of this algorithm is identical to
that of the Laplace case, it simply takes more operations (roughly twice as many) to
calculate the kernels and perform the required matrix operations.


We have assessed the accuracy of the algorithm for various domains, discretization
parameters, and Fourier modes.  The boundary conditions used are determined by
placing point charges inside the domains, and we evaluate the solution at random
points placed on a sphere that encompasses the boundary. In the combined field BIE
(\ref{eq:extDirBieHelmComb}), we set the parameter $\nu = k$.


\begin{table}[!ht]
\begin{center}
\begin{minipage}{0.99\linewidth} \begin{center}
 \begin{tabular}{| c | c | c | c | c | c |}
 \hline
 \multicolumn{1}{|c|}{$N_{\rm P}$} & \multicolumn{5}{c|}{$2N+1$} \\
 \hline
- &25&51&101&201&401    \\
   \hline
5&4.2306E-04&4.0187E-04&4.0185E-04&4.0185E-04&4.0185E-04\\
10&3.4791E-06&2.3738E-06&2.3759E-06&2.3760E-06&2.376E-06\\
20&7.8645E-06&4.5707E-09&7.1732E-11&7.1730E-11&7.173E-11\\
40&1.3908E-05&1.5980E-08&2.9812E-13&3.1500E-13&3.148E-13\\
80&1.0164E-05&5.4190E-09&4.9276E-13&4.8895E-13&-\\
  \hline
\end{tabular}
\end{center} \end{minipage} \\
\vspace{2mm}
\caption{Relative error in external Helmholtz problem for the domain in Figure \ref{fig:domains}(a). The domain is 1 wavelength in length (the major axis).  }
\label{tbl:helmerr1}
\end{center}
\end{table}


\begin{table}[!ht]
\begin{center}
\begin{minipage}{0.99\linewidth} \begin{center}
 \begin{tabular}{| c | c | c | c | c | c |}
 \hline
 \multicolumn{1}{|c|}{$N_{\rm P}$} & \multicolumn{5}{c|}{$2N+1$} \\
 \hline
- &25&51&101&201&401    \\
   \hline
5&2.0516E+00&3.4665E+00&4.1762E+00&4.4320E+00&4.4951E+00\\
10&2.3847E-01&2.4301E-01&2.4310E-01&2.4312E-01&2.4313E-01\\
20&1.7792E-02&8.8147E-06&8.8075E-06&8.8081E-06&8.8082E-06\\
40&1.7054E-02&5.4967E-08&3.7994E-10&3.7999E-10&3.7998E-10\\
80&1.7302E-02&1.6689E-08&1.8777E-11&1.8782E-11&-\\
  \hline
\end{tabular}
\end{center} \end{minipage} \\
\vspace{2mm}
\caption{Relative error in external Helmholtz problem for the domain in Figure \ref{fig:domains}(a). The domain is 25 wavelengths in length (the major axis).}
\label{tbl:helmerr2}
\end{center}
\end{table}

First, we consider the ellipsoidal domain given in Figure \ref{fig:domains}(a).  The
major axis of this ellipse has a diameter of 2, and its minor axes have a diameter
of 1/2. Table \ref{tbl:helmerr1} lists the accuracy achieved. We achieve 9 digits of
accuracy in this problem with 20 panels and 51 Fourier modes. We have not padded the
two sequences in the convolution procedure described in Section \ref{sec:kerExt},
but as more modes are used the tails of the sequence rapidly approach zero,
increasing the accuracy of the convolution algorithm utilized to calculate the
kernels. Table \ref{tbl:helmerr2} displays the same data, but with the wavenumber
increased so that there are 25 wavelengths along the length of the ellipsoid. We see
a minor decrease in accuracy as we would expect, but it only takes 40 panels  and 51
Fourier modes to achieve 8 digits of accuracy.


\begin{table}[!ht]
\begin{center}
\begin{minipage}{0.99\linewidth} \begin{center}
 \begin{tabular}{| c | c | c | c | c | c |}
 \hline
 \multicolumn{1}{|c|}{$N_{\rm P}$} & \multicolumn{5}{c|}{$2N+1$} \\
 \hline
- &25&51&101&201&401    \\
   \hline
5   & 2.6874E+00 & 2.5719E+00 & 2.5826E+00 & 2.2374E+00 & 1.9047E+00 \\
10 & 1.6351E+00 & 4.2972E+00 & 4.6017E+00 & 1.0380E+00 & 1.0328E+00 \\
20 & 4.3666E+00 & 2.6963E-03 & 2.6966E-03 & 2.6967E-03 & 2.6967E-03 \\
40 & 4.4498E+00 & 9.1729E-08 & 4.1650E-08 & 4.1661E-08 & 4.1664E-08 \\
80 & 4.3692E+00 & 7.3799E-08 & 1.4811E-10 & 1.4812E-10 & - \\
  \hline
\end{tabular}
\end{center} \end{minipage} \\
\vspace{2mm}
\caption{Relative error in external Helmholtz problem for the domain in Figure \ref{fig:domains}(b). The domain is 10 wavelengths in length (the major axis).}
\label{tbl:helmerr3}
\end{center}
\end{table}

\begin{table}[!ht]
\begin{center}
\begin{minipage}{0.99\linewidth} \begin{center}
 \begin{tabular}{| c | c | c | c | c | c |}
 \hline
 \multicolumn{1}{|c|}{$N_{\rm P}$} & \multicolumn{5}{c|}{$2N+1$} \\
 \hline
- &25&51&101&201&401    \\
   \hline
5   & 1.1140E+01 & 3.2400E+01 & 3.3885E+01 & 4.4322E+01 & 1.6151E+01 \\
10 & 1.9626E+01 & 7.7739E+01 & 6.2931E+01 & 3.6335E-01 & 3.6344E-01 \\
20 & 2.6155E+01 & 4.9485E+01 & 2.5796E+01 & 4.8239E-05 & 4.8245E-05 \\
40 & 3.3316E+01 & 5.0645E+01 & 2.4330E+01 & 1.3841E-09 & 1.3846E-09 \\
80 & 1.6966E+01 & 6.0163E+01 & 2.4354E+01 & 1.6510E-10 & - \\
  \hline
\end{tabular}
\end{center} \end{minipage} \\
\vspace{2mm}
\caption{Relative error in external Helmholtz problem for the domain in Figure \ref{fig:domains}(c). The domain is 10 wavelengths in length (the major axis).}
\label{tbl:helmerr4}
\end{center}
\end{table}

We now consider the more complex domains given in Figures \ref{fig:domains}(b) and
\ref{fig:domains}(c). They are a wavy shaped block and a starfish shaped block with the
outer diameter of size roughly $1.5$ and $1.0$, respectively. The accuracy for
various values of $N_{P}$ and $N$ are given in Table \ref{tbl:helmerr3} and
\ref{tbl:helmerr4}. We achieve 8 digits of accuracy with $40$ panels and $51$
Fourier modes for the wavy block and 9 digits of accuracy with $40$ panels and $201$
Fourier modes for the starfish block.

\lsp


\noindent
\textbf{\textit{Acknowledgements:}} The work reported was supported by
NSF grants DMS0748488 and DMS0941476.

\bibliography{main_bib,refs}

\providecommand{\bysame}{\leavevmode\hbox to3em{\hrulefill}\thinspace}
\providecommand{\MR}{\relax\ifhmode\unskip\space\fi MR }
\providecommand{\MRhref}[2]{%
  \href{http://www.ams.org/mathscinet-getitem?mr=#1}{#2}
}
\providecommand{\href}[2]{#2}
\begin{thebibliography}{10}

\bibitem{Abramowitz:65a}
M.~Abramowitz and I.A. Stegun, \emph{Handbook of mathematical functions with
  formulas, graphs, and mathematical tables}, Dover, New York, 1965.

\bibitem{alpert:99a}
Bradley~K. Alpert, \emph{Hybrid gauss-trapezoidal quadrature rules}, SIAM J.
  Sci. Comput. \textbf{20} (1999), 1551--1584.

\bibitem{Atkinson:97a}
K.~Atkinson, \emph{The numerical solution of integral equations of the second
  kind}, Cambridge University Press, Cambridge, 1997.

\bibitem{Bakr:85a}
A.A. Bakr, \emph{The boundary integral equation method in axisymmetric stress
  analysis problems}, Springer-Verlag, Berlin, 1985.

\bibitem{2011_bremer}
J.~Bremer, \emph{A fast direct solver for the integral equations of scattering
  theory on planar curves with corners}, Journal of Computational Physics
  (2011), no.~0, --.

\bibitem{Briggs:95a}
B.~Briggs and V.E. Henson, \emph{The {DFT}: An owner's manual for the discrete
  fourier transform}, {SIAM}, Philadelphia, 1995.

\bibitem{Bruno:01a}
O.P. Bruno and L.A. Kunyansky, \emph{A fast, high-order algorithm for the
  solution of surface scattering problems: basic implementation, test, and
  applications}, J. Comput. Phys. \textbf{169} (2001), 80--110.

\bibitem{Cohl:99a}
H.S. Cohl and J.E. Tohline, \emph{A compact cylindrical green's function
  expansion for the solution of potential problems}, Astrophys. J. \textbf{527}
  (1999), 86--101.

\bibitem{Fleming:04a}
J.L. Fleming, A.W. Wood, and W.D.~Wood Jr., \emph{Locally corrected nystr\"{o}m
  method for em scattering by bodies of revolution}, J. Comput. Phys.
  \textbf{196} (2004), 41--52.

\bibitem{Gil:07a}
A.~Gil, J.~Segura, and N.M. Temme, \emph{Numerical methods for special
  functions}, SIAM, Philadelphia, 2007.

\bibitem{Guenther:88a}
R.B. Guenther and J.W. Lee, \emph{Partial differential equations of
  mathematical physics and integral equations}, Dover, New York, 1988.

\bibitem{Gupta:1979a}
A.K. Gupta, \emph{The boundary integral equation method for potential problems
  involving axisymmetric geometry and arbitrary boundary conditions}, Master's
  thesis, University of Kentucky, 1979.

\bibitem{2011_hao_martinsson_quadrature}
S.~Hao, P.G. Martinsson, and P.~Young, \emph{High-order accurate nystrom
  discretization of integral equations with weakly singular kernels on smooth
  curves in the plane}, 2011, ar{X}iv.org report \#1112.6262.

\bibitem{Kapur:97a}
S.~Kapur and V.~Rokhlin, \emph{High-order corrected trapezoidal quadrature
  rules for singular functions}, SIAM J. Numer. Anal. \textbf{34} (1997),
  1331--1356.

\bibitem{2001_rokhlin_kolm}
P.~Kolm and V.~Rokhlin, \emph{Numerical quadratures for singular and
  hypersingular integrals}, Comput. Math. Appl. \textbf{41} (2001), 327--352.

\bibitem{Kress:83a}
R.~Kress and W.T. Spassov, \emph{On the condition of boundary integral
  operators for the exterior dirichlet problem for the {H}elmholtz equation},
  Numer. Math. \textbf{42} (1983), 77--95.

\bibitem{Kuijpers:97a}
A.H. Kuijpers, G.~Verbeek, and J.W. Verheij, \emph{An improved acoustic fourier
  boundary element method formulation using fast fourier transform
  integration}, J. Acoust. Soc. Am. \textbf{102} (1997), 1394--1401.

\bibitem{Martinsson:04a}
P.G. Martinsson and V.~Rokhlin, \emph{A fast direct solver for boundary
  integral equations in two dimensions}, J. Comput. Phys. \textbf{205} (2004),
  1--23.

\bibitem{Provatidis:98a}
C.~Provatidis, \emph{A boundary element method for axisymmetric potential
  problems with non-axisymmetric boundary conditions using fast fourier
  transform}, Engrg. Comput. \textbf{15} (1998), 428--449.

\bibitem{Rizzo:79a}
F.J. Rizzo and D.J. Shippy, \emph{A boundary integral approach to potential and
  elasticity problems for axisymmetric bodies with arbitrary boundary
  conditions}, Mech. Res. Commun. \textbf{6} (1979), 99--103.

\bibitem{Shippy:80a}
D.J. Shippy, F.J. Rizzo, and A.K. Gupta, \emph{Boundary-integral solution of
  potential problems involving axisymmetric bodies and nonsymmetric boundary
  conditions}, Developments in Theoretical and Applied Mechanics (J.E.
  Stoneking, ed.), 1980, pp.~189--206.

\bibitem{Soenarko:93a}
B.~Soenarko, \emph{A boundary element formuluation for radiation of acoustic
  waves from axisymmetric bodies with arbitrary boundary conditions}, J.
  Acoust. Soc. Am. \textbf{93} (1993), 631--639.

\bibitem{Tsinopoulos:99a}
S.V. Tsinopoulos, J.P. Agnantiaris, and D.~Polyzos, \emph{An advanced boundary
  element/fast fourier transform axisymmetric formulation for acoustic
  radiation and wave scattering problems}, J. Acoust. Soc. Am. \textbf{105}
  (1999), 1517--1526.

\bibitem{Wang:97a}
W.~Wang, N.~Atalla, and J.~Nicolas, \emph{A boundary integral approach for
  accoustic radiation of axisymmetric bodies with arbitrary boundary conditions
  valid for all wave numbers}, J. Acoust. Soc. Am. \textbf{101} (1997),
  1468--1478.

\end{thebibliography}
\bibliographystyle{amsplain}


\end{document}